\documentclass[12pt]{amsart}

\usepackage[paperheight=11in,paperwidth=8.5in,left=1in,top=1.25in,right=1in,bottom=1.25in,]{geometry}
\usepackage[pdfstartpage={1},pdfstartview={FitH},pdftitle={},pdfauthor={},pdfsubject={},pdfcreator={},pdfproducer={},pdfkeywords={},pagebackref=false,colorlinks=true,linkcolor=black,citecolor=black,filecolor=black,urlcolor=black,]{hyperref}

\usepackage{amsfonts}
\usepackage{amsmath}
\allowdisplaybreaks[4]
\usepackage{amssymb}
\usepackage{amsthm}
\usepackage{bbm}
\usepackage[justification=centering]{caption}
\usepackage[capitalize,nameinlink,noabbrev,nosort]{cleveref}
\usepackage{enumerate}
\usepackage{enumitem}
\usepackage{float}
\restylefloat{figure}
\usepackage{graphicx}
\usepackage{mathrsfs}
\usepackage[framemethod=tikz,ntheorem,xcolor]{mdframed}
\usepackage{parcolumns}
\usepackage{tabularx}
\usepackage{tikz}
\usepackage{tikz-cd}
\usetikzlibrary{arrows,calc,intersections,matrix,positioning,through}
\usepackage[all]{xy}
\usepackage[aligntableaux=center]{ytableau}

\numberwithin{equation}{section}

\newtheorem{thm}[equation]{Theorem}
\AfterEndEnvironment{thm}{\noindent\ignorespaces}
\newtheorem{cor}[equation]{Corollary}
\AfterEndEnvironment{cor}{\noindent\ignorespaces}
\newtheorem{lem}[equation]{Lemma}
\AfterEndEnvironment{lem}{\noindent\ignorespaces}
\newtheorem{prop}[equation]{Proposition}
\AfterEndEnvironment{prop}{\noindent\ignorespaces}

\theoremstyle{definition}
\newtheorem{defn}[equation]{Definition}
\AfterEndEnvironment{defn}{\noindent\ignorespaces}
\newtheorem*{pf_no_qed}{Proof}
\newenvironment{pf}[1][]{\begin{pf_no_qed}[#1]\pushQED{\qed}}{\popQED\end{pf_no_qed}}
\AfterEndEnvironment{pf}{\noindent\ignorespaces}
\newtheorem{eg_no_qed}[equation]{Example}
\newenvironment{eg}[1][]{\begin{eg_no_qed}[#1]\pushQED{\qed}}{\popQED\end{eg_no_qed}}
\AfterEndEnvironment{eg}{\noindent\ignorespaces}
\newtheorem{rmk}[equation]{Remark}
\AfterEndEnvironment{rmk}{\noindent\ignorespaces}

\theoremstyle{remark}
\newtheorem*{claim}{Claim}
\AfterEndEnvironment{claim}{\noindent\ignorespaces}
\newtheorem*{claimpf_no_qed}{Proof of Claim}
\newenvironment{claimpf}[1][]{\begin{claimpf_no_qed}[#1]\pushQED{\qed}}{\popQED\end{claimpf_no_qed}}
\AfterEndEnvironment{claimpf}{\noindent\ignorespaces}

\DeclareMathOperator{\alt}{alt}

\newcommand{\cd}{\cdots\!}

\DeclareMathOperator{\Gr}{Gr}

\newcommand{\rf}[1]{\hyperref[#1]{(\ref*{#1})}}
\DeclareMathOperator{\sgn}{sgn}
\DeclareMathOperator{\sign}{sign}

\DeclareMathOperator{\spn}{span}

\DeclareMathOperator{\var}{var}

\title{Sign variation, the Grassmannian, and total positivity}
\date{August 10th, 2016}

\author{Steven N. Karp}
\address{Department of Mathematics, University of California, Berkeley}
\email{\href{mailto:skarp@berkeley.edu}{skarp@berkeley.edu}}
\urladdr{\url{http://math.berkeley.edu/~skarp/}}

\thanks{This work was supported by a Chateaubriand fellowship, an NSERC postgraduate scholarship, and NSF grant DMS-1049513.}

\begin{document}

\begin{abstract}
The {\itshape totally nonnegative Grassmannian} is the set of $k$-dimensional subspaces $V$ of $\mathbb{R}^n$ whose nonzero Pl\"ucker coordinates all have the same sign. Gantmakher and Krein (1950) and Schoenberg and Whitney (1951) independently showed that $V$ is totally nonnegative iff every vector in $V$, when viewed as a sequence of $n$ numbers and ignoring any zeros, changes sign at most $k-1$ times. We generalize this result from the totally nonnegative Grassmannian to the entire Grassmannian, showing that if $V$ is {\itshape generic} (i.e.\ has no zero Pl\"ucker coordinates), then the vectors in $V$ change sign at most $m$ times iff certain sequences of Pl\"ucker coordinates of $V$ change sign at most $m-k+1$ times. We also give an algorithm which, given a non-generic $V$ whose vectors change sign at most $m$ times, perturbs $V$ into a generic subspace whose vectors also change sign at most $m$ times. We deduce that among all $V$ whose vectors change sign at most $m$ times, the generic subspaces are dense. These results generalize to oriented matroids. As an application of our results, we characterize when a generalized {\itshape amplituhedron} construction, in the sense of Arkani-Hamed and Trnka (2013), is well defined. We also give two ways of obtaining the {\itshape positroid cell} of each $V$ in the totally nonnegative Grassmannian from the sign patterns of vectors in $V$.
\end{abstract}

\maketitle

\section{Introduction and main results}\label{sec_introduction}

\noindent The {\itshape (real) Grassmannian $\Gr_{k,n}$} is the set of $k$-dimensional subspaces of $\mathbb{R}^n$. Given $V\in\Gr_{k,n}$, take a $k\times n$ matrix $X$ whose rows span $V$; then for $k$-subsets $I\subseteq\{1, \cd, n\}$, we let $\Delta_I(V)$ be the $k\times k$ minor of $X$ restricted to the columns in $I$, called a {\itshape Pl\"{u}cker coordinate}. (The $\Delta_I(V)$ depend on our choice of $X$ only up to a global constant.) If all nonzero $\Delta_I(V)$ have the same sign, then $V$ is called {\itshape totally nonnegative}, and if in addition no $\Delta_I(V)$ equals zero, then $V$ is called {\itshape totally positive}. For example, the span $V$ of $(1,0,0,-1)$ and $(-1,2,1,3)$ is a totally nonnegative element of $\Gr_{2,4}$, but $V$ is not totally positive since $\Delta_{\{2,3\}}(V) = 0$.

The set $\Gr_{k,n}^{\ge 0}$ of totally nonnegative $V\in\Gr_{k,n}$, called the {\itshape totally nonnegative Grassmannian}, has become a hot topic in algebraic combinatorics in the past two decades. The general algebraic study of total positivity for split reductive connected algebraic groups $G$ over $\mathbb{R}$, and partial flag varieties $G/P$, was initiated by Lusztig \cite{lusztig}, of which $\Gr_{k,n}^{\ge 0}$ corresponds to the special case $G/P = \Gr_{k,n}$. Of particular interest is the stratification of $\Gr_{k,n}^{\ge 0}$ according to whether each $\Delta_I$ is zero or nonzero. This stratification is a cell decomposition, which was conjectured by Lusztig \cite{lusztig} and proved by Rietsch \cite{rietsch} (for the general case $G/P$), and later understood combinatorially by Postnikov \cite{postnikov}.

This general theory traces its origin to the study of {\itshape totally positive matrices} in the 1930's, in the context of oscillation theory in analysis. Here positivity conditions on matrices can imply special oscillation and spectral properties. A well-known result of this kind is due to Gantmakher and Krein \cite{gantmakher_krein_1937}, which states that if an $n\times n$ matrix $X$ is {\itshape totally positive} (i.e.\ all $\binom{2n}{n}$ minors of $X$ are positive), then the $n$ eigenvalues of $X$ are distinct positive reals. Gantmakher and Krein \cite{gantmakher_krein_1950} also gave a characterization of (what would later be called) the totally nonnegative and totally positive Grassmannians in terms of sign variation. To state their result, we introduce some notation. For $v\in\mathbb{R}^n$, let $\var(v)$ be the number of times $v$ (viewed as a sequence of $n$ numbers, ignoring any zeros) changes sign, and let
$$
\overline{\var}(v) := \max\{\var(w) : \text{$w\in\mathbb{R}^n$ such that $w_i = v_i$ for all $1 \le i \le n$ with $v_i\neq 0$}\}.
$$
(We use the convention $\var(0) := -1$.) For example, if $v := (1, -1, 0, -2)\in\mathbb{R}^4$, then $\var(v) = 1$ and $\overline{\var}(v) = 3$.
\begin{thm}[Chapter V, Theorems 3 and 1 of \cite{gantmakher_krein_1950}]\label{gantmakher_krein}~\\
(i) $V\in\Gr_{k,n}$ is totally nonnegative iff $\var(v)\le k-1$ for all $v\in V$. \\
(ii) $V\in\Gr_{k,n}$ is totally positive iff $\overline{\var}(v)\le k-1$ for all $v\in V\setminus\{0\}$.
\end{thm}
(Part (i) above was proved independently by Schoenberg and Whitney \cite{schoenberg_whitney}.) For example, the two vectors $(1,0,0,-1)$ and $(-1,2,1,3)$ each change sign exactly once, and we can check that any vector in their span $V$ changes sign at most once, which is equivalent to $V$ being totally nonnegative. On the other hand, $\overline{\var}((1,0,0,-1)) = 3$, so $V$ is not totally positive. Every element of $\Gr_{k,n}$ has a vector which changes sign at least $k-1$ times (put a $k\times n$ matrix whose rows span $V$ into reduced row echelon form, and take the alternating sum of the rows), so the totally nonnegative elements are those whose vectors change sign as few times as possible.

The results of the paper are organized as follows. In \cref{sec_sign_changes}, we generalize \cref{gantmakher_krein} from the totally nonnegative Grassmannian to the entire Grassmannian, by giving a criterion for when $\var(v)\le m$ for all $v\in V$, or when $\overline{\var}(v)\le m$ for all $v\in V\setminus\{0\}$, in terms of the Pl\"{u}cker coordinates of $V$. (\cref{gantmakher_krein} is the case $m = k-1$.) As an application of our results, in \cref{sec_amplituhedron} we examine the construction of {\itshape amplituhedra} introduced by Arkani-Hamed and Trnka \cite{arkani-hamed_trnka}. In \cref{sec_positroids}, we show how to use the sign patterns of vectors in a totally nonnegative $V$ to determine the cell of $V$ in the cell decomposition of $\Gr_{k,n}^{\ge 0}$.

We briefly mention here that all of our results hold more generally for {\itshape oriented matroids}, and we prove them in this context. In this section we state our results in terms of the Grassmannian, so as to make them as accessible as possible. We introduce oriented matroids and the basic results about them we will need in \cref{sec_oriented_matroids}. See \cref{oriented_matroids_remark} at the end of this section for further comments about oriented matroids. \\

\noindent We now describe our main results. We let $[n] := \{1, 2, \cd, n\}$, and denote by $\binom{[n]}{r}$ the set of $r$-subsets of $[n]$.
\begin{thm}\label{intro_minor_criterion}
Suppose that $V\in\Gr_{k,n}$, and $m\ge k-1$. \\
(i) If $\var(v)\le m$ for all $v\in V$, then $\var((\Delta_{I\cup\{i\}}(V))_{i\in [n]\setminus I})\le m-k+1$ for all $I\in\binom{[n]}{k-1}$. \\
(ii) We have $\overline{\var}(v)\le m$ for all $v\in V\setminus\{0\}$ iff $\overline{\var}((\Delta_{I\cup\{i\}}(V))_{i\in [n]\setminus I})\le m-k+1$ for all $I\in\binom{[n]}{k-1}$ such that $\Delta_{I\cup\{i\}}(V)\neq 0$ for some $i\in [n]$.
\end{thm}
(See \cref{chirotope_criterion}.) If we take $m := k-1$, then we recover \cref{gantmakher_krein}; see \cref{gantmakher_krein_for_oriented_matroids} for the details.
\begin{eg}\label{minor_criterion_example}
Let $V\in\Gr_{2,4}$ be the row span of the matrix $\begin{bmatrix}1 & 0 & -2 & 3 \\ 0 & 2 & 1 & 4\end{bmatrix}$, so $k := 2$. Then by \cref{intro_minor_criterion}(ii), the fact that $\overline{\var}(v)\le 2 =: m$ for all $v\in V\setminus\{0\}$ is equivalent to the fact that the $4$ sequences
\begin{align*}
& (\Delta_{\{1,2\}}(V), \Delta_{\{1,3\}}(V), \Delta_{\{1,4\}}(V)) = (2, 1, 4), \\
& (\Delta_{\{1,2\}}(V), \Delta_{\{2,3\}}(V), \Delta_{\{2,4\}}(V)) = (2, 4, -6), \\
& (\Delta_{\{1,3\}}(V), \Delta_{\{2,3\}}(V), \Delta_{\{3,4\}}(V)) = (1, 4, -11), \\
& (\Delta_{\{1,4\}}(V), \Delta_{\{2,4\}}(V), \Delta_{\{3,4\}}(V)) = (4, -6, -11)
\end{align*}
each change sign at most $m-k+1=1$ time.
\end{eg}
We say that $V\in\Gr_{k,n}$ is {\itshape generic} if all Pl\"{u}cker coordinates of $V$ are nonzero. If $V$ is generic, then (ii) above implies that the converse of (i) holds. The converse of (i) does not hold in general (see \cref{non_generic_counterexample}); however, if $V\in\Gr_{k,n}$ is not generic and $\var(v)\le m$ for all $v\in V$, then we show how to perturb $V$ into a generic $V'\in\Gr_{k,n}$ while maintaining the property $\var(v)\le m$ for all $v\in V'$. Working backwards, we can then apply \cref{intro_minor_criterion}(i) to $V'$ in order to test whether $\var(v)\le m$ for all $v\in V$. The precise statement is as follows.
\begin{thm}\label{intro_perturbation_criterion}
Given $V\in\Gr_{k,n}$, we can perturb $V$ into a generic $V'\in\Gr_{k,n}$ such that $\max_{v\in V}\var(v)$ = $\max_{v\in V'}\var(v)$. In particular, for $m\ge k-1$ we have $\var(v)\le m$ for all $v\in V$ iff $\var((\Delta_{I\cup\{i\}}(V'))_{i\in [n]\setminus I})\le m-k+1$ for all $I\in\binom{[n]}{k-1}$.

Thus in $\{V\in\Gr_{k,n} : \var(v)\le m \text{ for all } v\in V\}$, the generic elements are dense.
\end{thm}
(See \cref{perturbation_criterion} and \cref{generic_elements_are_dense}.) In the special case $m = k-1$, we recover the result of Postnikov (Section 17 of \cite{postnikov}) that the totally positive Grassmannian is dense in the totally nonnegative Grassmannian.

\cref{perturbation_criterion} in fact gives an algorithm for perturbing $V$ into a generic $V'$. It involves taking a $k\times n$ matrix $X$ whose rows span $V$, and repeatedly adding a very small multiple of a column of $X$ to an adjacent column (and taking the row span of the resulting matrix). We show that repeating the sequence $1\to_+ 2, 2\to_+ 3, \cd, (n-1)\to_+ n, n\to_+ (n-1), (n-1)\to_+ (n-2), \cd, 2\to_+ 1$ of adjacent-column perturbations $k$ times in order from left to right is sufficient to obtain a generic $V'$, where $i\to_+ j$ denotes adding a very small positive multiple of column $i$ to column $j$. We give several other sequences of adjacent-column perturbations which work; see \cref{perturbation_criterion}.

We use these results to study {\itshape amplituhedra}, introduced by Arkani-Hamed and Trnka \cite{arkani-hamed_trnka} to help calculate scattering amplitudes in theoretical physics. They consider the map $\Gr_{k,n}^{\ge 0}\to\Gr_{k,r}$ on the totally nonnegative Grassmannian induced by a given linear map $Z:\mathbb{R}^n\to\mathbb{R}^r$. Note that this map is not necessarily well defined, since $Z$ may send a $k$-dimensional subspace to a subspace of lesser dimension. In part to preclude this possibility, Arkani-Hamed and Trnka require that $k\le r$ and $Z$ has positive $r\times r$ minors (when viewed as an $r\times n$ matrix), and call the image of the map $\Gr_{k,n}^{\ge 0}\to\Gr_{k,r}$ a {\itshape (tree) amplituhedron}. Lam \cite{lam} showed more generally that the map $\Gr_{k,n}^{\ge 0}\to\Gr_{k,r}$ is well defined if the row span of $Z$ (regarded as an $r\times n$ matrix) has a totally positive $k$-dimensional subspace, in which case he calls the image a {\itshape (full) Grassmann polytope}. (When $k=1$, Grassmann polytopes are precisely polytopes in projective space.) We use sign variation to give a necessary and sufficient condition for the map $\Gr_{k,n}^{\ge 0}\to\Gr_{k,r}$ induced by $Z$ to be well defined, and in particular recover the sufficient conditions of Arkani-Hamed and Trnka, and Lam.
\begin{thm}\label{intro_amplituhedron_map}
Suppose that $k,n,r\in\mathbb{N}$ with $n\ge k,r$, and that $Z:\mathbb{R}^n\to\mathbb{R}^{r}$ is a linear map, which we also regard as an $r\times n$ matrix. Let $d$ be the rank of $Z$ and $W\in\Gr_{d,n}$ the row span of $Z$, so that $W^\perp = \ker(Z)\in\Gr_{n-d,n}$. The following are equivalent: \\
(i) the map $\Gr_{k,n}^{\ge 0}\to\Gr_{k,r}$ induced by $Z$ is well defined, i.e.\ $\dim(Z(V)) = k$ for all $V\in\Gr_{k,n}^{\ge 0}$; \\
(ii) $\var(v)\ge k$ for all nonzero $v\in\ker(Z)$; and \\
(iii) $\overline{\var}((\Delta_{I\setminus\{i\}}(W))_{i\in I})\le d-k$ for all $I\in\binom{[n]}{d+1}$ such that $W|_I$ has dimension $d$.
\end{thm}
(See \cref{amplituhedron_map}.) We remark that the equivalence of (ii) and (iii) above is equivalent to \cref{intro_minor_criterion}(ii). \\

\noindent We now describe our results about the cell decomposition of $\Gr_{k,n}^{\ge 0}$. Given $V\in\Gr_{k,n}$, we define the {\itshape matroid $M(V)$} of $V$ as the set of $I\in\binom{[n]}{k}$ such that $\Delta_I(V)$ is nonzero. If $V$ is totally nonnegative, we also call $M(V)$ a {\itshape positroid}. The stratification of $\Gr_{k,n}^{\ge 0}$ by positroids (i.e.\ its partition into equivalence classes, where $V\sim W$ iff $M(V) = M(W)$) is a cell decomposition \cite{rietsch, postnikov}.

How can we determine the matroid of $V\in\Gr_{k,n}$ from the sign patterns of vectors in $V$? Given $I\subseteq [n]$ and a sign vector $\omega\in\{+,-\}^I$, we say that {\itshape $V$ realizes $\omega$} if there exists a vector in $V$ whose restriction to $I$ has signs given by $\omega$. For example, if $(2,3,-2,-1)\in V$, then $V$ realizes $(+,-,-)$ on $\{1,3,4\}$. Note that $V$ realizes $\omega$ iff $V$ realizes $-\omega$. It is not difficult to show that for all $I\in\binom{[n]}{k}$, we have $I\in M(V)$ iff $V$ realizes all $2^k$ sign vectors in $\{+,-\}^I$. Furthermore, in order to determine whether $I$ is in $M(V)$ from which sign vectors $V$ realizes in $\{+,-\}^I$, we potentially have to check all $2^{k-1}$ pairs of sign vectors (each sign vector and its negation), since given any $\omega\in\{+,-\}^I$ (and assuming $n > k$), there exists $V\in\Gr_{k,n}$ which realizes all $2^k$ sign vectors in $\{+,-\}^I$ except for $\pm\omega$. (See \cref{check_all_patterns}.) However, in the case that $V$ is totally nonnegative, we show that we need only check $k$ particular sign vectors in $\{+,-\}^I$ to verify that $\Delta_I(V)\neq 0$.
\begin{thm}\label{intro_basis_criterion}
For $V\in\Gr_{k,n}^{\ge 0}$ and $I\in\binom{[n]}{k}$, we have $I\in M(V)$ iff $V$ realizes the $2k$ (or $k$ up to sign) sign vectors in $\{+,-\}^I$ which alternate in sign between every pair of consecutive components, with at most one exceptional pair. 
\end{thm}
(See \cref{basis_criterion}.) For example, if $k = 5$, these $2k$ sign vectors are $(+,-,+,-,+)$, $(+,+,-,+,-)$, $(+,-,-,+,-)$, $(+,-,+,+,-)$, $(+,-,+,-,-)$, and their negations.
\begin{eg}\label{basis_criterion_example}
Let $V\in\Gr_{3,5}^{\ge 0}$ be the row span of the matrix $\begin{bmatrix}2 & 1 & 0 & 0 & 3 \\ 0 & 0 & 1 & 0 & 0 \\ 0 & 0 & 0 & 1 & 1\end{bmatrix}$. \cref{intro_basis_criterion} implies that for all $I\in\binom{[5]}{3}$, we have $\Delta_I(V)\neq 0$ iff $V$ realizes the $3$ sign vectors $(+,-,+)$, $(+,+,-)$, $(+,-,-)$ on $I$. For $I = \{1,3,5\}$, the vectors $(2,1,-1,0,3)$, $(2,1,1,-4,-1)$, $(2,1,-1,-4,-1)\in V$ realize the sign vectors $(+,-,+)$, $(+,+,-)$, $(+,-,-)$ on $I$, so $\Delta_{\{1,3,5\}}(V)\neq 0$. (We do not need to check that $(+,+,+)$, the remaining sign vector in $\{+,-\}^I$ up to sign, is realized.) For $I = \{1,4,5\}$, the vectors $(2,1,0,-1,2), (2,1,0,-4,-1)\in V$ realize the sign vectors $(+,-,+), (+,-,-)$ on $I$, but no vector in $V$ realizes the sign vector $(+,+,-)$ on $I$, so $\Delta_{\{1,4,5\}}(V) = 0$.
\end{eg}

We now describe another way to recover the positroid of $V\in\Gr_{k,n}^{\ge 0}$ from the sign patterns of vectors in $V$. We begin by showing how to obtain the {\itshape Schubert cell} of $V$, which is labeled by the lexicographic minimum of $M(V)$. To state this result, we introduce some notation. For $v\in\mathbb{R}^n$ and $I\subseteq [n]$, we say that $v$ {\itshape strictly alternates in sign on $I$} if $v|_I$ has no zero components, and alternates in sign between consecutive components. Let $A(V)$ denote the set of $I\in\binom{[n]}{k}$ such that some vector in $V$ strictly alternates in sign on $I$. Note that if $I\in M(V)$ then $V|_I = \mathbb{R}^I$, so $M(V)\subseteq A(V)$. We also define the {\itshape Gale partial order} $\le_\textnormal{Gale}$ on $\binom{[n]}{k}$ by $\{i_1 < \cdots < i_k\}\le_\textnormal{Gale}\{j_1 < \cdots < j_k\}$ iff $i_1 \le j_1, i_2 \le j_2, \cd, i_k \le j_k$.
\begin{thm}\label{intro_schubert_cell_criterion}
For $V\in\Gr_{k,n}^{\ge 0}$, the lexicographic minimum of $M(V)$ equals the Gale minimum of $A(V)$.
\end{thm}
(See \cref{schubert_cell_criterion}.) We remark that the lexicographic minimum of $M(V)$ is also the Gale minimum of $M(V)$, but $A(V)$ does not necessarily equal $M(V)$ (see \cref{alternating_set_can_be_bigger}). We also note that if $V\in\Gr_{k,n}$ is not totally nonnegative, then $A(V)$ does not necessarily have a Gale minimum (see \cref{gale_minimality_example}).
\begin{eg}\label{alternating_set_can_be_bigger}
Let $V\in\Gr_{3,5}^{\ge 0}$ be the row span of the matrix $\begin{bmatrix}2 & 1 & 0 & 0 & 3 \\ 0 & 0 & 1 & 0 & 0 \\ 0 & 0 & 0 & 1 & 1\end{bmatrix}$, as in \cref{basis_criterion_example}. \cref{intro_schubert_cell_criterion} implies that the lexicographic minimum $\{1,3,4\}$ of $M(V)$ equals the Gale minimum of
$$
A(V) = \{\{1,3,4\},\{1,3,5\},\{1,4,5\},\{2,3,4\},\{2,3,5\},\{2,4,5\},\{3,4,5\}\}.
$$
Note that $\{2,4,5\}\in A(V)\setminus M(V)$.
\end{eg}
By the cyclic symmetry of the totally nonnegative Grassmannian, we can then use \cref{intro_schubert_cell_criterion} to recover the {\itshape Grassmann necklace} of $V\in\Gr_{k,n}^{\ge 0}$ (see \cref{grassmann_necklace_criterion}), which in turn determines the positroid of $V$ by a result of Postnikov (Theorem 17.1 of \cite{postnikov}).
\begin{rmk}
We can easily reinterpret results about upper bounds on $\var$ in terms of lower bounds on $\overline{\var}$, and upper bounds on $\overline{\var}$ in terms of lower bounds on $\var$, by the following two facts.
\begin{lem}\label{dual_translation}
(i) \cite{gantmakher_krein_1950} For $v\in\mathbb{R}^n\setminus\{0\}$, we have
$$
\var(v) + \overline{\var}(\alt(v)) = n-1,
$$
where $\alt(v) := (v_1, -v_2, v_3, -v_4, \cd, (-1)^{n-1}v_n)\in\mathbb{R}^n$. \\(ii) \cite{hilbert}\cite{hochster} Given $V\in\Gr_{k,n}$, let $V^\perp\in\Gr_{n-k,n}$ be the orthogonal complement of $V$. Then $V$ and $\alt(V^\perp)$ have the same Pl\"{u}cker coordinates:
$$
\Delta_I(V) = \Delta_{[n]\setminus I}(\alt(V^\perp)) \qquad \text{ for all } I\in\binom{[n]}{k}.
$$
\end{lem}
(Part (i) is stated without proof as equation (67) in Chapter II of \cite{gantmakher_krein_1950}; see equation (5.1) of \cite{ando} for a proof. The earliest satement of part (ii) we found in the literature is at the beginning of Section 7 of \cite{hochster}. Hochster does not give a proof, and says that this result ``was basically known to Hilbert.'' The idea is that if $[I_k | A]$ is a $k\times n$ matrix whose rows span $V\in\Gr_{k,n}$, where $A$ is a $k\times (n-k)$ matrix, then $V^\perp$ is the row span of the matrix $[A^T | -I_{n-k}]$. This idea appears implicitly in equation (14) of \cite{hilbert}, and more explicitly in Theorem 2.2.8 of \cite{oxley} and Proposition 3.1(i) of \cite{musiker_reiner}. I thank a referee for pointing out the reference \cite{musiker_reiner}.) For example, we get the following dual formulation of Gantmakher and Krein's result (\cref{gantmakher_krein}).
\begin{cor}[[Chapter V, Theorems 7 and 6 of \cite{gantmakher_krein_1950}]\label{dual_gantmakher_krein} ~\\
(i) $V\in\Gr_{k,n}$ is totally nonnegative iff $\overline{\var}(v) \ge k$ for all $v\in V\setminus\{0\}$. \\
(ii) $V\in\Gr_{k,n}$ is totally positive iff $\var(v)\ge k$ for all $v\in V\setminus\{0\}$.
\end{cor}
\end{rmk}

\begin{rmk}\label{oriented_matroids_remark}
The natural framework in which to consider sign patterns of vectors in $V$, and signs of the Pl\"{u}cker coordinates of $V$, is that of {\itshape oriented matroids}. Our results hold, and are proven, in this context, and are more general because while every subspace gives rise to an oriented matroid, not every oriented matroid comes from a subspace. (The totally nonnegative Grassmannian is a special case; the analogue of a totally nonnegative subspace is a {\itshape positively oriented matroid}, and Ardila, Rinc\'{o}n, and Williams \cite{ardila_rincon_williams} recently showed that each positively oriented matroid comes from a totally nonnegative subspace. Hence there is no added generality gained here in passing from the Grassmannian to oriented matroids.)

Those already familiar with oriented matroids can use the following dictionary to reinterpret the results stated in this section:
\begin{center}
\begin{tabular}{c|c}
subspaces & oriented matroids \\ \hline\hline
sign vectors of vectors in $V$ & {\itshape covectors} of $\mathcal{M}(V)$ \\ \hline
$\Delta(V)$, up to sign & the {\itshape chirotope} $\chi_{\mathcal{M}(V)}$ \\ \hline
$V$ is generic & $\mathcal{M}(V)$ is {\itshape uniform} \\ \hline
$V'$ is a perturbation of $V$ & there is a weak map from $\mathcal{M}(V')$ to $\mathcal{M}(V)$ \\ \hline
the closure of $S\subseteq\Gr_{k,n}$ & images of weak maps from $\mathcal{M}(V')$, over $V'\in S$ \\ \hline
the orthogonal complement $V^\perp$ of $V$ & the dual of $\mathcal{M}(V)$ \\ \hline
\end{tabular}
\end{center}
We also generalize to oriented matroids the operation of adding a very small multiple of a column (of a $k\times n$ matrix whose rows span $V\in\Gr_{k,n}$) to an adjacent column; see \cref{defn_perturbation}.

For those unfamiliar with oriented matroids, we give an introduction in \cref{sec_oriented_matroids}, biased toward the tools we need. For a thorough introduction to oriented matroids, see the book \cite{bjorner_las_vergnas_sturmfels_white_ziegler}.
\end{rmk}

\noindent{\bfseries Acknowledgements.} I thank Lauren Williams, my advisor, for many helpful conversations and suggestions. I also thank the referees for their valuable feedback. I am grateful to Sylvie Corteel and the Laboratoire d'Informatique Algorithmique: Fondements et Applications at Universit\'{e} Paris Diderot for hosting me while I conducted part of this work.

\section{Introduction to oriented matroids}\label{sec_oriented_matroids}

\noindent In this section, we introduce oriented matroids, and much of the notation and tools that we use in our proofs. A comprehensive account of the theory of oriented matroids, and our reference throughout, is the book by Bj\"{o}rner, Las Vergnas, Sturmfels, White, and Ziegler \cite{bjorner_las_vergnas_sturmfels_white_ziegler}. We begin by describing oriented matroids coming from subspaces of $\mathbb{R}^n$ (i.e.\ {\itshape realizable} oriented matroids), which will serve as motivation for the exposition to follow.

For $\alpha\in\mathbb{R}$ we define
$$
\sign(\alpha) := \begin{cases}
0, & \text{if $\alpha = 0$} \\
+, & \text{if $\alpha > 0$} \\
-, & \text{if $\alpha < 0$}
\end{cases},
$$
and for $x\in\mathbb{R}^E$ we define the sign vector $\sign(x)\in\{0,+,-\}^E$ by $\sign(x)_e := \sign(x_e)$ for $e\in E$. We will sometimes use $1$ and $-1$ in place of $+$ and $-$. For example, $\sign(5,0,-1,2) = (+, 0, -, +) = (1, 0, -1, 1)$. Given a sign vector $X\in\{0,+,-\}^E$, the {\itshape support} of $X$ is the subset $\underline{X} := \{e\in E : X_e \neq 0\}$ of $E$. We can think of $X$ as giving a sign to each element of $\underline{X}$ (some authors call $X$ a {\itshape signed subset}). We also define $-X\in\{0,+,-\}^E$ by $(-X)_e := -X_e$ for $e\in E$. For example, $X = (+, 0, -, +)\in\{0,+,-\}^4$ has support $\{1,3,4\}$, and $-X = (-, 0, +, -)$.
\begin{defn}[realizable oriented matroids; 1.2 of \cite{bjorner_las_vergnas_sturmfels_white_ziegler}]\label{defn_realizable_oriented_matroids}
Let $E$ be a finite set and $V$ a $k$-dimensional subspace of $\mathbb{R}^E$. The {\itshape (realizable) oriented matroid $\mathcal{M}(V)$} associated to $V$ is uniquely determined by $E$ (the {\itshape ground set} of $\mathcal{M}(V)$) and any one of the following three objects: \\
$\bullet$ the set $\mathcal{V}^* := \{\sign(v) : v\in V\}$, called the {\itshape covectors} of $\mathcal{M}(V)$; or \\
$\bullet$ the set $\mathcal{C}^* := \{X\in \mathcal{V}^* : X \text{ has minimal nonempty support}\}$, called the {\itshape cocircuits} of $\mathcal{M}(V)$; or \\
$\bullet$ the function $\chi : E^k \to \{0,+,-\}$ (up to multiplication by $\pm 1$), called the {\itshape chirotope} of $\mathcal{M}(V)$, where $\chi(i_1, \cd, i_k) := \sign(\det([x^{(i_1)} | \cdots | x^{(i_k)}]))$ ($i_1, \cd, i_k\in E$) for some fixed $k\times E$ matrix $[x^{(i)} : i\in E]$ whose rows span $V$. \\
The {\itshape rank} of $\mathcal{M}(V)$ is $k$.
\end{defn}

\begin{eg}\label{eg_realizable_matroid}
Let $V\in\Gr_{2,3}$ be the row span of the matrix $\begin{bmatrix}0 & -1 & 1 \\ 3 & 0 & 2\end{bmatrix}$. Then $\mathcal{M}(V)$ is an oriented matroid of rank $k := 2$ with ground set $E := \{1,2,3\}$. Note that $(+, +, -)$ is a covector of $\mathcal{M}(V)$, because it is the sign vector of e.g.\ $(3, 3, -1)\in V$. The covectors of $\mathcal{M}(V)$ are
\begin{gather*}
(0, 0, 0), (0, +, -), (0, -, +), (+, 0, +), (-, 0, -), (+, +, 0), (-, -, 0), \\
(+, +, -), (-, +, -), (+, -, +), (-, -, +), (+, +, +), (-, -, -).
\end{gather*}
The cocircuits of $\mathcal{M}(V)$ are the covectors with minimal nonempty support, i.e.\ $$
(0, +, -), (0, -, +), (+, 0, +), (-, 0, -), (+, +, 0), (-, -, 0).
$$
The chirotope $\chi$ of $\mathcal{M}(V)$ is given (up to sign) by
\begin{align*}
& \chi(1,2) = \sign(\Delta_{\{1,2\}}(V)) = \sign(3) = +, \\
& \chi(1,3) = \sign(\Delta_{\{1,3\}}(V)) = \sign(-3) = -, \\
& \chi(2,3) = \sign(\Delta_{\{2,3\}}(V)) = \sign(-2) = -, 
\end{align*}
and by the fact that $\chi$ is {\itshape alternating}, i.e.\ swapping two arguments multiplies the result by $-1$. The fact that the Pl\"{u}cker coordinates $\Delta_I(V)$ are defined only up to multiplication by a global nonzero constant explains why the chirotope is defined only up to sign.
\end{eg}
\begin{defn}[oriented matroid, cocircuit axioms; 3.2.1 of \cite{bjorner_las_vergnas_sturmfels_white_ziegler}]\label{defn_cocircuits}
An {\itshape oriented matroid $\mathcal{M}$} is an ordered pair $\mathcal{M} = (E,\mathcal{C}^*)$, where $E$ is a finite set and $\mathcal{C}^*\subseteq 2^{\{0, +, -\}^E}$ satisfies the following four axioms: \\
(C0) every sign vector in $\mathcal{C}^*$ has nonempty support; \\
(C1) $\mathcal{C}^* = -\mathcal{C}^*$; \\
(C2) if $X,Y\in\mathcal{C}^*$ with $\underline{X}\subseteq\underline{Y}$, then $X = \pm Y$; \\
(C3) if $X,Y\in\mathcal{C}^*$ and $a\in E$ such that $X\neq -Y$ and $X_a = -Y_a \neq 0$, then there exists $Z\in\mathcal{C}^*$ such that $Z_a = 0$, and $Z_b = X_b$ or $Z_b = Y_b$ for all $b\in\underline{Z}$.

The set $E$ is called the {\itshape ground set} of $\mathcal{M}$, and the sign vectors in $\mathcal{C}^*$ are called the {\itshape cocircuits} of $\mathcal{M}$.
\end{defn}
We denote the cocircuits of $\mathcal{M}$ by $\mathcal{C}^*(\mathcal{M})$. (The superscript $*$ is present to indicate that cocircuits are {\itshape circuits} of the {\itshape dual} of $\mathcal{M}$.) We remark that not all oriented matroids are realizable; see 1.5.1 of \cite{bjorner_las_vergnas_sturmfels_white_ziegler} for an example of a non-realizable oriented matroid.
\begin{eg}\label{eg_cocircuits}
The sign vectors
$$
(0, +, -), (0, -, +), (+, 0, +), (-, 0, -), (+, -, 0), (-, +, 0)
$$
are {\itshape not} the cocircuits of an oriented matroid, because e.g.\ (C3) above fails when we take $X = (0, +, -)$, $Y = (+, 0, +)$, $a = 3$.
\end{eg}

For sign vectors $X,Y\in\{0,+,-\}^E$, define the {\itshape composition $X\circ Y$} as the sign vector in $\{0,+,-\}^E$ given by
$$
(X\circ Y)_e := \begin{cases}
X_e, & \text{if $X_e\neq 0$} \\
Y_e, & \text{if $X_e = 0$}
\end{cases} \qquad \text{ for } e\in E.
$$
We can think of $X\circ Y$ as being formed by starting with $X$ and recording $Y$ in the empty slots of $X$, or by starting with $Y$ and overwriting $X$ on top. In general, $X\circ Y\neq Y\circ X$; if the composition of sign vectors $X^{(1)}, \cd, X^{(r)}$ of $E$ does not depend on the order of composition, we say that $X^{(1)}, \cd, X^{(r)}$ are {\itshape conformal}. For example, $(+,0,-)$ and $(0,+,-)$ are conformal.

A {\itshape covector} of an oriented matroid $\mathcal{M}$ is a composition of some (finite number of) cocircuits of $\mathcal{M}$. (We include the empty composition, which is the zero sign vector.) We let $\mathcal{V}^*(\mathcal{M})$ denote the set of covectors of $\mathcal{M}$. Note that by (C2) of \cref{defn_cocircuits}, we can recover the cocircuits of $\mathcal{M}$ as the covectors with minimal nonempty support. A key property of covectors is the following conformality property.
\begin{prop}[conformality for covectors; 3.7.2 of \cite{bjorner_las_vergnas_sturmfels_white_ziegler}]\label{conformality_for_covectors}
Suppose that $X$ is a covector of the oriented matroid $\mathcal{M}$. Then $X = C^{(1)} \circ \cdots \circ C^{(r)}$ for some conformal cocircuits $C^{(1)}, \cd, C^{(r)}$ of $\mathcal{M}$.
\end{prop}
There are axioms which characterize when a set of sign vectors in $\{0,+,-\}^E$ is the set of covectors of an oriented matroid; see 3.7.5 of \cite{bjorner_las_vergnas_sturmfels_white_ziegler}.
\begin{defn}[basis, rank; pp124, 115 of \cite{bjorner_las_vergnas_sturmfels_white_ziegler}]\label{defn_basis_rank}
Let $\mathcal{M}$ be an oriented matroid with ground set $E$. A {\itshape basis of $\mathcal{M}$} is a minimal $B\subseteq E$ such that $B\cap\underline{C}\neq\emptyset$ for every cocircuit $C$ of $\mathcal{M}$. All bases of $\mathcal{M}$ have the same size $k\ge 0$, called the {\itshape rank} of $\mathcal{M}$.
\end{defn}
$\mathcal{M}$ determines a unique orientation on its bases (up to a global sign).
\begin{defn}[chirotope; 3.5.1, 3.5.2 of \cite{bjorner_las_vergnas_sturmfels_white_ziegler}]\label{defn_chirotope}
Suppose that $\mathcal{M}$ is an oriented matroid of rank $k$ with ground set $E$. Then there exists a function $\chi_\mathcal{M} : E^k\to\{0,+,-\}$ (called the {\itshape chirotope} of $\mathcal{M}$), unique up to sign, satisfying the following properties: \\
(i) $\chi_\mathcal{M}$ is {\itshape alternating}, i.e.\ $\chi_\mathcal{M}(i_{\sigma(1)}, \cd, i_{\sigma(k)}) = \sgn(\sigma)\chi_\mathcal{M}(i_1, \cd, i_k)$ for $i_1, \cd, i_k\in E$ and $\sigma\in\mathfrak{S}_k$; \\
(ii) $\chi_\mathcal{M}(i_1, \cd, i_k) = 0$ if $\{i_1, \cd, i_k\}\subseteq E$ is not a basis of $\mathcal{M}$; and \\
(iii) if $\{a, i_1, \cd, i_{k-1}\}, \{b, i_1, \cd, i_{k-1}\}\subseteq E$ are bases of $\mathcal{M}$ and $C$ is a cocircuit of $\mathcal{M}$ with $i_1, \cd, i_{k-1}\notin\underline{C}$, then $\chi_\mathcal{M}(a, i_1, \cd, i_{k-1}) = C_aC_b\chi_\mathcal{M}(b, i_1, \cd, i_{k-1})$.
\end{defn}
$\mathcal{M}$ is uniquely determined by $\chi_\mathcal{M}$ up to sign (3.5.2 of \cite{bjorner_las_vergnas_sturmfels_white_ziegler}). For the axioms characterizing chirotopes of oriented matroids, see 3.5.3, 3.5.4 of \cite{bjorner_las_vergnas_sturmfels_white_ziegler}. If $E$ is totally ordered (for example, $E = [n]$ ordered by $1 < 2 < \cdots < n$), we let $\chi_\mathcal{M}(I)$ denote $\chi_\mathcal{M}(i_1, \cd, i_k)$ for $I\in\binom{E}{k}$ ($I = \{i_1, \cd, i_k\}$, $i_1 < \cdots < i_k$), and set $\chi_\mathcal{M}(J) := 0$ for $J\subseteq E$ with $|J| < k$. In this case $\chi_\mathcal{M}$ gives an {\itshape orientation} (either $+$ or $-$) to each basis of $\mathcal{M}$.

The relation (iii) above between $\chi_\mathcal{M}$ and the cocircuits of $\mathcal{M}$ is called the {\itshape pivoting property}; we state it in the following useful form.
\begin{prop}[pivoting property; 3.5.1, 3.5.2 of \cite{bjorner_las_vergnas_sturmfels_white_ziegler}]\label{pivoting_property}
Suppose that $\mathcal{M}$ is an oriented matroid of rank $k$ with a totally ordered ground set $E$, $I\in\binom{E}{k-1}$, and $a, b\in E$. If $I\cup\{a\}$ and $I\cup\{b\}$ are bases of $\mathcal{M}$, then there exists a cocircuit $C$ of $\mathcal{M}$ with $I\cap\underline{C} = \emptyset$ (unique up to sign), whence $a,b\in\underline{C}$, and
\begin{align}\label{pivoting_property_equation}
\chi_\mathcal{M}(I\cup\{a\}) = (-1)^{|\{i\in I : i \textnormal{ is strictly between $a$ and $b$}\}|}C_aC_b\chi_\mathcal{M}(I\cup\{b\}).
\end{align}
Conversely, if there exists a cocircuit $C$ of $\mathcal{M}$ with $I\cap\underline{C} = \emptyset$ and $b\in\underline{C}$, then \rf{pivoting_property_equation} holds.
\end{prop}
Only the first part of \cref{pivoting_property} is proved in \cite{bjorner_las_vergnas_sturmfels_white_ziegler}, so we prove the converse.
\begin{pf}[of converse]
Let $C$ be a cocircuit of $\mathcal{M}$ with $I\cap\underline{C} = \emptyset$ and $b\in\underline{C}$. First suppose that $I\cup\{b\}$ is not a basis of $\mathcal{M}$; we must show that $I\cup\{a\}$ is also not a basis. By \cref{defn_basis_rank} there exists a cocircuit $D$ of $\mathcal{M}$ with $(I\cup\{b\})\cap\underline{D} = \emptyset$. If $a\notin\underline{C}$ or $a\notin\underline{D}$, then we immediately get that $I\cup\{a\}$ is not a basis. Otherwise we have $a\in\underline{C}\cup\underline{D}$, and $C\neq\pm D$ since $b\in\underline{C}\setminus\underline{D}$. Hence we may apply (C3) of \cref{defn_cocircuits} to obtain a cocircuit of $\mathcal{M}$ whose support is contained in $(\underline{C}\cup\underline{D})\setminus\{a\}\subseteq E\setminus (I\cup\{a\})$, whence $I\cup\{a\}$ is not a basis of $\mathcal{M}$. Similarly, if $a\in\underline{C}$ and $I\cup\{a\}$ is not a basis of $\mathcal{M}$, then $I\cup\{b\}$ is not a basis, giving \rf{pivoting_property_equation}. Also, if $I\cup\{a\}$ and $I\cup\{b\}$ are both bases of $\mathcal{M}$, then \rf{pivoting_property_equation} follows from the first part of this result. The remaining case is when $I\cup\{b\}$ is a basis of $\mathcal{M}$, $I\cup\{a\}$ is not a basis, and $a\notin\underline{C}$, whence both sides of \rf{pivoting_property_equation} are zero.
\end{pf}

Now we introduce restriction of oriented matroids; for a realizable oriented matroid $\mathcal{M}(V)$, this corresponds to restricting $V$ to a subset of the canonical coordinates.
\begin{defn}[restriction; 3.7.11, 3.4.9, pp133-134 of \cite{bjorner_las_vergnas_sturmfels_white_ziegler}]\label{defn_restriction}
Let $\mathcal{M}$ be an oriented matroid with ground set $E$, and $F\subseteq E$. The {\itshape restriction of $\mathcal{M}$ to $F$}, denoted by $\mathcal{M}|_F$ or $\mathcal{M}\setminus G$ (where $G := E\setminus F$), is the oriented matroid with ground set $F$ and covectors $\{X|_F : X\in\mathcal{V}^*(\mathcal{M})\}$. The bases of $\mathcal{M}|_F$ are the maximal elements of $\{B\cap F : B\text{ is a basis of }\mathcal{M}\}$. The chirotope $\chi_{\mathcal{M}|_F}$ is given as follows. Let $k,l$ be the ranks of $\mathcal{M},\mathcal{M}|_F$, respectively, and take $i_1, \cd, i_{k-l}\in E\setminus F$ such that $F\cup\{i_1, \cd, i_{k-l}\}$ contains a basis of $\mathcal{M}$. Then
$$
\chi_{\mathcal{M}|_F}(j_1, \cd, j_l) = \chi_\mathcal{M}(j_1, \cd, j_l, i_1, \cd, i_{k-l}) \qquad \text{ for } j_1, \cd, j_l\in F.
$$
If $V$ is a subspace of $\mathbb{R}^E$, then $\mathcal{M}(V)|_F = \mathcal{M}(V|_F)$.
\end{defn}

We conclude by describing a partial order on oriented matroids with a fixed ground set. Geometrically, for point configurations, moving up in the partial order corresponds to moving the points of the configuration into more general position. (A configuration of $n$ points in $\mathbb{R}^k$ gives rise to a subspace of $\mathbb{R}^n$, and hence an oriented matroid with ground set $[n]$, by writing the points as the columns of a $k\times n$ matrix and taking the row span of this matrix.) We use the partial order on sign vectors given by $X\le Y$ iff $Y = X\circ Y$ ($X,Y\in\{0,+,-\}^E$), i.e.\ $X_e = Y_e$ for all $e\in E$ such that $X_e\neq 0$. This also defines a partial order on chirotopes, regarded as sign vectors in $\{0,+,-\}^{E^k}$.
\begin{defn}[partial order on oriented matroids; 7.7.5 of \cite{bjorner_las_vergnas_sturmfels_white_ziegler}]\label{defn_poset}
Let $\mathcal{M}$, $\mathcal{N}$ be oriented matroids with ground set $E$. We say that $\mathcal{M}\le\mathcal{N}$ if for every covector $X$ of $\mathcal{M}$, there exists a covector $Y$ of $\mathcal{N}$ with $X\le Y$. Then $\le$ is a partial order on oriented matroids with ground set $E$. If $\mathcal{M}$ and $\mathcal{N}$ have the same rank, then $\mathcal{M}\le\mathcal{N}$ iff $\chi_\mathcal{M}\le\pm\chi_\mathcal{N}$.
\end{defn}
The standard terminology for $\mathcal{M}\le\mathcal{N}$ is that there is a {\itshape weak map} from $\mathcal{N}$ to $\mathcal{M}$.

\section{Relating sign changes of covectors and the chirotope}\label{sec_sign_changes}

\noindent Recall that given a sign vector $X\in\{0,+,-\}^E$ over a totally ordered set $E$, the number of sign changes of $X$ (ignoring any zeros) is denoted by $\var(X)$, and $\overline{\var}(X) := \max_{Y\ge X}\var(Y)$. The goal of this section is to give, for any oriented matroid $\mathcal{M}$ with a totally ordered ground set, a criterion for when $\var(X)\le m$ for all covectors $X$ of $\mathcal{M}$, or when $\overline{\var}(X)\le m$ for all nonzero covectors $X$ of $\mathcal{M}$, in terms of the chirotope of $\mathcal{M}$. \cref{chirotope_criterion} provides such a criterion in the latter case, as well as in the former case if $\mathcal{M}$ is {\itshape uniform}, i.e.\ every $k$-subset of its ground set is a basis (where $k$ is the rank of $\mathcal{M}$). (Hence $V\in\Gr_{k,n}$ is generic iff $\mathcal{M}(V)$ is uniform.) For non-uniform $\mathcal{M}$, we then show (\cref{perturbation_criterion}) how to perturb $\mathcal{M}$ into a generic uniform matroid $\mathcal{N}$ so that we may apply the criterion in \cref{chirotope_criterion} to determine whether $\var(X)\le m$ for all covectors $X$ of $\mathcal{M}$.

We remark that while $\var$ is weakly increasing (i.e.\ $\var(X)\le\var(Y)$ if $X\le Y$), $\overline{\var}$ is weakly decreasing, which helps explain why $\var$ and $\overline{\var}$ require such different treatments.
\begin{thm}\label{chirotope_criterion}
Suppose that $\mathcal{M}$ is an oriented matroid of rank $k$ with ground set $[n]$, and $m \ge k-1$. \\
(i) If $\var(X)\le m$ for all $X\in\mathcal{V}^*(\mathcal{M})$, then $\var((\chi_\mathcal{M}(I\cup\{i\}))_{i\in [n]\setminus I})\le m-k+1$ for all $I\in\binom{[n]}{k-1}$. \\
(ii) We have $\overline{\var}(X)\le m$ for all $X\in\mathcal{V}^*(\mathcal{M})\setminus\{0\}$ iff $\overline{\var}((\chi_\mathcal{M}(I\cup\{i\}))_{i\in [n]\setminus I})\le m-k+1$ for all $I\in\binom{[n]}{k-1}$ such that $I\cup\{i\}$ is a basis of $\mathcal{M}$ for some $i\in [n]$.
\end{thm}
For an example using this theorem, see \cref{minor_criterion_example}. Note that (ii) above implies that if $\mathcal{M}$ is uniform, then the converse of (i) holds. However, the converse of (i) does not hold in general, as shown in \cref{non_generic_counterexample}. \cref{exclude_zero_sequences_counterexample} shows that the condition ``$I\cup\{i\}$ is a basis of $\mathcal{M}$ for some $i\in [n]$'' (equivalently, that the sequence $(\chi_\mathcal{M}(I\cup\{i\}))_{i\in [n]\setminus I}$ is nonzero) in (ii) is necessary. Also note that there is no loss of generality in the assumption $m\ge k-1$, because there exists a covector of $\mathcal{M}$ which changes sign at least $k-1$ times; in fact, if $B\in\binom{[n]}{k}$ is any basis of $\mathcal{M}$, then there exists a covector of $\mathcal{M}$ which strictly alternates in sign on $B$. (This follows from \cref{defn_restriction}: $\mathcal{M}|_B$ is the uniform oriented matroid of rank $k$ with ground set $B$, and so $\mathcal{V}^*(\mathcal{M})|_B = \mathcal{V}^*(\mathcal{M}|_B) = \{0,+,-\}^B$.)
\begin{eg}\label{non_generic_counterexample}
Let $V\in\Gr_{2,4}$ be the row span of the matrix $\begin{bmatrix}1 & 0 & 1 & 0 \\ 0 & 1 & 0 & 1\end{bmatrix}$, so $k := 2$. Note that the $4$ sequences of Pl\"{u}cker coordinates
\begin{align*}
&(\Delta_{\{1,2\}}(V) , \Delta_{\{1,3\}}(V) , \Delta_{\{1,4\}}(V)) = (1 , 0 , 1), \\
&(\Delta_{\{1,2\}}(V) , \Delta_{\{2,3\}}(V) , \Delta_{\{2,4\}}(V)) = (1 , -1 , 0), \\
&(\Delta_{\{1,3\}}(V) , \Delta_{\{2,3\}}(V) , \Delta_{\{3,4\}}(V)) = (0 , -1 , 1), \\
&(\Delta_{\{1,4\}}(V) , \Delta_{\{2,4\}}(V) , \Delta_{\{3,4\}}(V)) = (1 , 0 , 1)
\end{align*}
each change sign at most $m-k+1=1$ time (where we take $m := 2$), but the vector $(1,-1,1,-1)\in V$ changes sign $3$ times. Hence the converse to \cref{chirotope_criterion}(i) does not hold. However, if we were forced to pick a sign for, say, $\Delta_{\{1,3\}}(V)$, then either the first or third sequence above would change sign twice. This motivates the introduction of perturbations below.
\end{eg}

\begin{eg}\label{exclude_zero_sequences_counterexample}
Let $V\in\Gr_{3,5}$ be the row span of the matrix $\begin{bmatrix}1 & 1 & 0 & 0 & 0 \\ 0 & 0 & 1 & 0 & -1 \\ 0 & 0 & 0 & 1 & 1\end{bmatrix}$. Then $\mathcal{M}(V)$ satisfies the equivalent conditions of \cref{chirotope_criterion}(ii) with $m := 3$, i.e.\ $\overline{\var}(v) \le 3$ for all $v\in V\setminus\{0\}$, and $\overline{\var}((\Delta_{I\cup\{i\}}(V))_{i\in [5]\setminus I})\le 1$ for all $I\in\binom{[5]}{2}$ such that $\Delta_{I\cup\{i\}}(V)\neq 0$ for some $i\in [n]$. We cannot remove the condition ``$\Delta_{I\cup\{i\}}(V)\neq 0$ for some $i\in [n]$,'' because taking $J := \{1,2\}$ we have $\Delta_{J\cup\{i\}}(V) = 0$ for all $i\in [n]$, and so $\overline{\var}((\Delta_{J\cup\{i\}}(V))_{i\in [5]\setminus J}) = 2$.
\end{eg}
\begin{pf}[of \cref{chirotope_criterion}]
The idea is to use the fact that if $X\in\{0,+,-\}^n$ with $\var(X) = r$, then there exists $A\in\binom{[n]}{r+1}$ such that $X$ strictly alternates in sign on $A$. We restrict our attention to an appropriate choice of $A$, using \cref{defn_basis_rank} and the pivoting property (\cref{pivoting_property}) to relate cocircuits, bases, and the chirotope.

(i) Suppose that $I\in\binom{[n]}{k-1}$ such that $\var((\chi_\mathcal{M}(I\cup\{i\}))_{i\in [n]\setminus I})\ge m-k+2$. Take $A\in\binom{[n]\setminus I}{m-k+3}$ such that $(\chi_\mathcal{M}(I\cup\{i\}))_{i\in [n]\setminus I}$ strictly alternates in sign on $A$. Fix $a\in A$, and for the remainder of this proof, for $i,j\in [n]$ let $[i,j)$ denote the interval of integers from $i$ (inclusively) to $j$ (exclusively), i.e.\ $\{i, i+1, \cd, j-1\}$ if $j \ge i$ and $\{j+1, j+2, \cd, i\}$ if $j \le i$. By \cref{defn_basis_rank}, for $i\in I$ there exists a cocircuit $C^{(i)}$ of $\mathcal{M}$ with $((I\cup\{a\})\setminus\{i\})\cap\underline{C^{(i)}} = \emptyset$; since $I\cup\{a\}$ is a basis of $\mathcal{M}$ we have $i\in\underline{C^{(i)}}$, so we may assume that $C^{(i)}_i = (-1)^{|(I\cup A)\cap[a,i)|}$. Also let $D$ be a cocircuit of $\mathcal{M}$ with $I\cap\underline{D} = \emptyset$; since $a\in\underline{D}$ we may assume that $D_a = 1$. Then for $b\in [n]$, the pivoting property (\cref{pivoting_property}) gives
$$
\chi_\mathcal{M}(I\cup\{b\}) = (-1)^{|I\cap [a,b)|}D_aD_b\chi_\mathcal{M}(I\cup\{a\}).
$$
Because $(\chi_\mathcal{M}(I\cup\{i\}))_{i\in A}$ strictly alternates in sign, we have
$$
\chi_\mathcal{M}(I\cup\{b\}) = (-1)^{|A\cap [a,b)|}\chi_\mathcal{M}(I\cup\{a\}) \qquad \text{ for } b\in A,
$$
so $D_b = (-1)^{|(I\cup A)\cap [a,b)|}$. Now let $X$ be the covector $D\circ C^{(i_1)}\circ\cdots\circ C^{(i_{k-1})}$ of $\mathcal{M}$, where $I = \{i_1, \cd, i_{k-1}\}$. Then $X_i = (-1)^{|(I\cup A)\cap [a,i)|}$ for $i\in I\cup A$, so $X$ strictly alternates in sign on $I\cup A$, giving $\var(X)\ge m+1$.

(ii) ($\Rightarrow$): Suppose that $I\in\binom{[n]}{k-1}$ such that $I\cup\{i\}$ is a basis of $\mathcal{M}$ for some $i\in [n]$, and $\overline{\var}((\chi_\mathcal{M}(I\cup\{i\}))_{i\in [n]\setminus I})\ge m-k+2$. We proceed as in the proof of (i). Take $A\in\binom{[n]\setminus I}{m-k+3}$ such that $\overline{\var}((\chi_\mathcal{M}(I\cup\{i\}))_{i\in A}) = m-k+2$ and $I\cup\{a\}$ is a basis of $\mathcal{M}$ for some $a\in A$; fix such an $a\in A$. By \cref{defn_basis_rank} there exists a cocircuit $D$ of $\mathcal{M}$ with $I\cap\underline{D} = \emptyset$; since $a\in\underline{D}$ we may assume that $D_a = 1$. Then for $b\in [n]$, the pivoting property (\cref{pivoting_property}) gives
$$
\chi_\mathcal{M}(I\cup\{b\}) = (-1)^{|I\cap [a,b)|}D_aD_b\chi_\mathcal{M}(I\cup\{a\}).
$$
Because $\overline{\var}((\chi_\mathcal{M}(I\cup\{i\}))_{i\in A}) = m-k+2$, for $b\in A$ either $\chi_\mathcal{M}(I\cup\{b\}) = 0$ or
$$
\chi_\mathcal{M}(I\cup\{b\}) = (-1)^{|A\cap [a,b)|}\chi_\mathcal{M}(I\cup\{a\}),
$$
whence either $D_b = 0$ or $D_b = (-1)^{|(I\cup A)\cap [a,b)|}$. Hence $D|_{I\cup A} \le X$, where $X$ is the sign vector in $\{0,+,-\}^{I\cup A}$ with $X_a = 1$ which strictly alternates in sign. This gives $\overline{\var}(D) \ge \var(X) = m+1$.

($\Leftarrow$): Suppose that $\overline{\var}(X)\ge m+1$ for some nonzero covector $X$ of $\mathcal{M}$. By \cref{conformality_for_covectors} there exists a cocircuit $C$ of $\mathcal{M}$ with $C \le X$, whence $\overline{\var}(C)\ge m+1$. We consider two cases. First suppose that $|\underline{C}|\le n-m-1$. Take $a\in\underline{C}$, and note that by (C2) of \cref{defn_cocircuits}, $([n]\setminus\underline{C})\cup\{a\}$ has nonempty intersection with the support of every cocircuit of $\mathcal{M}$. Hence by \cref{defn_basis_rank}, some subset of $([n]\setminus\underline{C})\cup\{a\}$ is a basis of $\mathcal{M}$, which we may write as $I\cup\{a\}$ for some $I\in\binom{[n]\setminus\underline{C}}{k-1}$. Then $I\cup\{i\}$ is not a basis of $\mathcal{M}$ for $i\in [n]\setminus\underline{C}$, whence $(\chi_\mathcal{M}(I\cup\{i\}))_{i\in [n]\setminus I}$ has at least $m-k+2$ zero components. This gives $\overline{\var}((\chi_\mathcal{M}(I\cup\{i\}))_{i\in [n]\setminus I}) \ge m-k+2$, completing the proof.

Now suppose instead that $|\underline{C}|\ge n-m-1$. There exists $J\in\binom{[n]}{m+2}$ with $\overline{\var}(C|_J) = m+1$; take such a $J$ which minimizes $|J\cap\underline{C}|$. It follows that $[n]\setminus\underline{C}\subseteq J$. (Otherwise there exists $e\in [n]\setminus (J\cup\underline{C})$, whence letting $e'$ equal either $\min_{f\in J\cap\underline{C}, f > e}f$ or $\max_{f\in J\cap\underline{C}, f < e}f$, at least one of which exists because $|\underline{C}|\ge n-m-1$, we have $\overline{\var}(C|_{(J\setminus\{e'\})\cup\{e\}}) = m+1$, contradicting our choice of $J$.) Since $|\underline{C}|\ge n-m-1$, we may take $j\in J\cap\underline{C}$. Note that by (C2) of \cref{defn_cocircuits}, $([n]\setminus\underline{C})\cup\{j\}$ has nonempty intersection with the support of every cocircuit of $\mathcal{M}$. Hence by \cref{defn_basis_rank}, some subset of $([n]\setminus\underline{C})\cup\{j\}$ is a basis of $\mathcal{M}$, which we may write as $I\cup\{j\}$ for some $I\in\binom{[n]\setminus\underline{C}}{k-1}$. In particular, we have $I\subseteq J$.

By the pivoting property (\cref{pivoting_property}), we have
$$
\chi_\mathcal{M}(I\cup\{i\}) = (-1)^{|I\cap [j,i)|}C_iC_j\chi_\mathcal{M}(I\cup\{j\}) \qquad \text{ for } i\in [n].
$$
Also, since $C$ weakly alternates in sign on $J$, for $i\in J$ we have either $C_i = 0$ or $C_i = (-1)^{J\cap [j,i)}C_j$. Hence for $i\in J\setminus I$ we have either $\chi_\mathcal{M}(I\cup\{i\}) = 0$ or $\chi_\mathcal{M}(I\cup\{i\}) = (-1)^{(J\setminus I)\cap [j,i)}\chi_\mathcal{M}(I\cup\{j\})$, whence $(\chi_\mathcal{M}(I\cup\{i\}))_{i\in J\setminus I}$ weakly alternates in sign on $J\setminus I$, i.e.\ $\overline{\var}((\chi_\mathcal{M}(I\cup\{i\}))_{i\in J\setminus I}) = m-k+2$.
\end{pf}

We call an oriented matroid $\mathcal{M}$ with a totally ordered ground set {\itshape positively oriented} if every basis of $\mathcal{M}$ has the same orientation, and {\itshape alternating} if $\mathcal{M}$ is positively oriented and uniform. Hence $V\in\Gr_{k,n}$ is totally nonnegative iff $\mathcal{M}(V)$ is positively oriented, and $V$ is totally positive iff $\mathcal{M}(V)$ is alternating. We now obtain the generalization of Gantmakher and Krein's characterization (\cref{gantmakher_krein}) to oriented matroids, as a consequence of \cref{chirotope_criterion} in the special case $m := k-1$.

\begin{cor}\label{gantmakher_krein_for_oriented_matroids}
Suppose that $\mathcal{M}$ is an oriented matroid of rank $k$ with ground set $[n]$. \\
(i) $\mathcal{M}$ is positively oriented iff $\var(X)\le k-1$ for all $X\in\mathcal{V}^*(\mathcal{M})$. \\
(ii) $\mathcal{M}$ is alternating iff $\overline{\var}(X)\le k-1$ for all $X\in\mathcal{V}^*(\mathcal{M})\setminus\{0\}$.
\end{cor}
We remark that the forward directions of (i) and (ii) above follow from \cref{gantmakher_krein} and Ardila, Rinc\'{o}n, and Williams' result \cite{ardila_rincon_williams} that every positively oriented matroid is realizable. (The converses do not so follow, because we do not know {\itshape a priori} that an oriented matroid $\mathcal{M}$ satisfying $\var(X)\le k-1$ for all $X\in\mathcal{V}^*(\mathcal{M})$ is realizable.) Part (ii) above is implicit in the literature (cf.\ \cite{bland_las_vergnas}, 9.4 of \cite{bjorner_las_vergnas_sturmfels_white_ziegler}, and \cite{cordovil_duchet}), though we have not seen it explicitly stated and proven in this form.
\begin{pf}
(i) ($\Rightarrow$): Suppose that $\mathcal{M}$ is positively oriented, and let $\mathcal{N}$ be the uniform positively oriented matroid of rank $k$ with ground set $[n]$. Then by \cref{chirotope_criterion}(ii) with $m := k-1$, we have $\overline{\var}(Y)\le k-1$ for all $Y\in\mathcal{V}^*(\mathcal{N})\setminus\{0\}$. Now given any $X\in\mathcal{V}^*(\mathcal{M})$, since $\mathcal{M}\le\mathcal{N}$ there exists $Y\in\mathcal{V}^*(\mathcal{N})$ with $X\le Y$ (\cref{defn_poset}), whence $\var(X)\le\var(Y)\le k-1$.

($\Leftarrow$): Suppose that $\var(X)\le k-1$ for all $X\in\mathcal{V}^*(\mathcal{M})$. Then by \cref{chirotope_criterion}(i) with $m := k-1$, any two bases of $\mathcal{M}$ which have $k-1$ elements in common have the same orientation. Hence it will suffice to show that given any two bases $I$, $J$ of $\mathcal{M}$, there exist bases $I_0 := I, I_1, \cd, I_{r-1}, I_r := J$ of $\mathcal{M}$ such that $|I_{s-1}\cap I_s| \ge k-1$ for all $s\in [r]$. This follows from the {\itshape basis exchange} axiom for (oriented) matroids (p81 of \cite{bjorner_las_vergnas_sturmfels_white_ziegler}): if $A$ and $B$ are bases of an (oriented) matroid and $a\in A\setminus B$, then there exists $b\in B\setminus A$ such that $(A\setminus\{a\})\cup\{b\}$ is a basis.

(ii) The forward direction follows from \cref{chirotope_criterion}(ii) with $m := k-1$. For the converse, suppose that $\overline{\var}(X)\le k-1$ for all $X\in\mathcal{V}^*(\mathcal{M})\setminus\{0\}$. Then $\mathcal{M}$ is positively oriented by part (i) of this result. Also, if there exists $I\in\binom{[n]}{k}$ which is not a basis of $\mathcal{M}$, then by \cref{defn_basis_rank} there exists a cocircuit $C$ of $\mathcal{M}$ with $I\cap\underline{C} = \emptyset$, whence $\overline{\var}(C)\ge k$, a contradiction. Hence $\mathcal{M}$ is uniform.
\end{pf}

We have already observed that the converse to \cref{chirotope_criterion}(i) holds when $\mathcal{M}$ is a uniform oriented matroid, but not in general. Our goal in the remainder of the section is to prove a necessary and sufficient condition for having $\var(X)\le m$ for all $X\in\mathcal{V}^*(\mathcal{M})$. Namely, we give an algorithm for perturbing any oriented matroid $\mathcal{M}$ with a totally ordered ground set into a uniform $\mathcal{N}\ge\mathcal{M}$ of the same rank, such that $\max_{X\in\mathcal{V}^*(\mathcal{M})}\var(X) = \max_{Y\in\mathcal{V}^*(\mathcal{N})}\var(Y)$; we then apply \cref{chirotope_criterion} to $\mathcal{N}$ to determine $\max_{X\in\mathcal{V}^*(\mathcal{M})}\var(X)$ (\cref{perturbation_criterion}). In the case of realizable oriented matroids $\mathcal{M}(V)$ ($V\in\Gr_{k,n}$), this perturbation involves repeatedly adding a very small multiple of one column of a $k\times n$ matrix whose rows span $V$ to an adjacent column (and taking the row span of the resulting matrix). These perturbations generalize to all oriented matroids, as we explain below.

Let $\mathcal{M}$ be an oriented matroid with ground set $E$. A {\itshape single element extension of $\mathcal{M}$ at $a$} is an oriented matroid $\widetilde{\mathcal{M}}$ with ground set $E\sqcup\{a\}$ (where $\sqcup$ denotes disjoint union) and the same rank as $\mathcal{M}$, such that $\widetilde{\mathcal{M}}|_E = \mathcal{M}$. (Some authors allow $\widetilde{\mathcal{M}}$ to have rank greater than $\mathcal{M}$.) Las Vergnas \cite{las_vergnas_1978} studied single element extensions; we use his results as stated in \cite{bjorner_las_vergnas_sturmfels_white_ziegler}. For a sign vector $X\in\{0,+,-\}^E$ and $y\in\{0,+,-\}$, let $(X,y)_a\in\{0,+,-\}^{E\sqcup\{a\}}$ denote the sign vector whose restriction to $E$ is $X$ and whose $a$th component is $y$.

\begin{lem}[cocircuits of single element extensions; 7.1.4 of \cite{bjorner_las_vergnas_sturmfels_white_ziegler}]\label{single_element_cocircuits}
Suppose that the oriented matroid $\widetilde{\mathcal{M}}$ is the single element extension of $\mathcal{M}$ at $a$, where $\mathcal{M}$ has ground set $E$ and rank $k$. Then there exists a unique function $\sigma:\mathcal{C}^*(\mathcal{M})\to\{0,+,-\}$ such that $(C,\sigma(C))_a$ is a cocircuit of $\widetilde{\mathcal{M}}$ for all cocircuits $C$ of $\mathcal{M}$. We have $\mathcal{C}^*(\widetilde{\mathcal{M}}) = \mathcal{C}_1\sqcup\mathcal{C}_2$, where
\begin{align*}
& \mathcal{C}_1 := \{(C,\sigma(C))_a : C\in\mathcal{C}^*(\mathcal{M})\}, \\
& \mathcal{C}_2 := \left\{(C\circ D,0)_a : \begin{aligned}
& C,D\in\mathcal{C}^*(\mathcal{M}) \textnormal{ are conformal}, \sigma(C) = -\sigma(D) \neq 0, \textnormal{ and} \\
& |B\setminus(\underline{C}\cup\underline{D})| \ge k-2 \textnormal{ for some basis } B \textnormal{ of } \mathcal{M}
\end{aligned}\right\}.
\end{align*}
\end{lem}
In this case we say that $\widetilde{\mathcal{M}}$ is the single element extension of $\mathcal{M}$ at $a$ {\itshape by $\sigma$}. In general, not all functions $\sigma:\mathcal{C}^*(\mathcal{M})\to\{0,+,-\}$ give rise to single element extensions. However, for $e\in E$ the evaluation function $\phi_e:\mathcal{C}^*(\mathcal{M})\to\{0,+,-\}, C\mapsto C_e$ and its negation $-\phi_e$ are guaranteed to give single element extensions (7.1.8 of \cite{bjorner_las_vergnas_sturmfels_white_ziegler}). (Geometrically, the single element extension by $\phi_e$ duplicates coordinate $e$ at the new coordinate $a$.) Also, if the two functions $\sigma, \tau : \mathcal{C}^*(\mathcal{M})\to\{0,+,-\}$ each give rise to a single element extension of $\mathcal{M}$, then so does the composition $\sigma\circ\tau$ (7.2.2 of \cite{bjorner_las_vergnas_sturmfels_white_ziegler}), and the extension of $\mathcal{M}$ by $\sigma$ is less than or equal to (under \cref{defn_poset}) the extension of $\mathcal{M}$ by $\sigma\circ\tau$ (7.7.8 of \cite{bjorner_las_vergnas_sturmfels_white_ziegler}). (The composition $\sigma\circ\tau:\mathcal{C}^*(\mathcal{M})\to\{0,+,-\}$ is defined just as for covectors, by
$$
(\sigma\circ\tau)(C) := \begin{cases}
\sigma(C), & \text{if $\sigma(C)\neq 0$} \\
\tau(C), & \text{if $\sigma(C) = 0$}
\end{cases} \qquad \text{ for } C\in\mathcal{C}^*(\mathcal{M}).)
$$
\begin{defn}[$i\to_\epsilon j$-perturbation]\label{defn_perturbation}
Let $\mathcal{M}$ be an oriented matroid with ground set $E$, $i,j\in E$, and $\epsilon\in\{+,-\}$. If $i = j$ or $j$ is a coloop of $\mathcal{M}$, set $\mathcal{N} := \mathcal{M}$. (A {\itshape coloop} $c$ of an (oriented) matroid is an element of its ground set which is in every basis.) Otherwise, the restriction $\mathcal{M}\setminus\{j\}$ has the same rank as $\mathcal{M}$, so $\mathcal{M}$ is the single element extension of $\mathcal{M}\setminus\{j\}$ at $j$ by some $\sigma:\mathcal{C}^*(\mathcal{M}\setminus\{j\})\to\{0,+,-\}$. Let $\mathcal{N}$ be the single element extension of $\mathcal{M}\setminus\{j\}$ at $j$ by $\sigma\circ\epsilon\phi_i$, which is well defined and satisfies $\mathcal{N}\ge\mathcal{M}$ by the preceding discussion. We call $\mathcal{N}$ the {\itshape $i\to_\epsilon j$-perturbation of $\mathcal{M}$}.
\end{defn}

We now prove several properties of $i\to_\epsilon j$-perturbation.
\begin{lem}[chirotope of the $i\to_\epsilon j$-perturbation]\label{perturbation_chirotope}
Suppose that $\mathcal{M}$ is an oriented matroid of rank $k$ with a totally ordered ground set $E$, and $\mathcal{N}$ is the $i\to_\epsilon j$-perturbation of $\mathcal{M}$ (where $i,j\in E$ and $\epsilon\in\{+,-\}$). Then the chirotope of $\mathcal{N}$ is given by
$$
\chi_\mathcal{N}(I) = \begin{cases}
(-1)^{|I\cap (i,j)|}\epsilon\chi_\mathcal{M}((I\setminus\{j\})\cup\{i\}), & \textnormal{if $i\notin I$, $j\in I$, and $\chi_\mathcal{M}(I) = 0$} \\
\chi_\mathcal{M}(I), & \textnormal{otherwise}
\end{cases}
$$
for $I\in\binom{E}{k}$, where $(i,j)$ denotes the set of elements of $E$ strictly between $i$ and $j$.
\end{lem}

\begin{pf}
If $i=j$ or $j$ is a coloop of $\mathcal{M}$, then $\mathcal{N} = \mathcal{M}$ and the result is clear, so we may assume that $i\neq j$ and $j$ is not a coloop of $\mathcal{M}$. Let $I\in\binom{E}{k}$. If $j\notin I$, we have $\chi_\mathcal{N}(I) = \chi_\mathcal{M}(I)$ because $\mathcal{M}\setminus\{j\} = \mathcal{N}\setminus\{j\}$. Also, if $\chi_\mathcal{M}(I)\neq 0$, we have $\chi_\mathcal{N}(I) = \chi_\mathcal{M}(I)$ because $\mathcal{M}\le\mathcal{N}$. Hence we may assume that $j\in I$ and $\chi_\mathcal{M}(I) = 0$. Then by \cref{defn_basis_rank} there exists a cocircuit $C$ of $\mathcal{M}$ with $I\cap\underline{C} = \emptyset$. In particular, $C_j = 0$. If $C_i = 0$, then by \cref{single_element_cocircuits} $C$ is also a cocircuit of $\mathcal{N}$, whence by \cref{defn_basis_rank} both $I$ and $(I\setminus\{j\})\cup\{i\}$ are not bases of $\mathcal{M}$ or $\mathcal{N}$, giving $\chi_\mathcal{N}(I) = \chi_\mathcal{M}((I\setminus\{j\})\cup\{i\}) = \chi_\mathcal{M}(I) = 0$.

Suppose instead that $C_i \neq 0$. In particular $i\notin I$, so we must show that $\chi_\mathcal{N}(I) = (-1)^{|I\cap (i,j)|}\epsilon\chi_\mathcal{M}((I\setminus\{j\})\cup\{i\})$. By \cref{single_element_cocircuits}, we get a cocircuit $D$ of $\mathcal{N}$ such that $D_e = C_e$ for $e\in E\setminus\{j\}$, and either $D_j = \epsilon C_i$ (if $C\in\mathcal{C}_1$) or $D_j = 0$ (if $C\in\mathcal{C}_2$). Hence by the pivoting property (\cref{pivoting_property}), we have
$$
\chi_\mathcal{N}(I) = (-1)^{|I\cap (i,j)|}D_iD_j\chi_\mathcal{N}((I\setminus\{j\})\cup\{i\}),
$$
and $\chi_\mathcal{N}((I\setminus\{j\})\cup\{i\}) = \chi_\mathcal{M}((I\setminus\{j\})\cup\{i\})$ since $j\notin(I\setminus\{j\})\cup\{i\}$. If $C\in\mathcal{C}_1$, then $D_iD_j = \epsilon$, giving $\chi_\mathcal{N}(I) = (-1)^{|I\cap (i,j)|}\epsilon\chi_\mathcal{M}((I\setminus\{j\})\cup\{i\})$. Now suppose that $C\in\mathcal{C}_2$. Then $D_j = 0$, giving $\chi_\mathcal{N}(I) = 0$; we must show that $\chi_\mathcal{M}((I\setminus\{j\})\cup\{i\}) = 0$. Since $C\in\mathcal{C}_2$, we can write $C = (X\circ Y,0)_j$ for some conformal cocircuits $X,Y$ of $\mathcal{M}\setminus\{j\}$ with $\sigma(X) = -\sigma(Y) \neq 0$, where $\mathcal{M}$ is the single element extension of $\mathcal{M}\setminus\{j\}$ by $\sigma$. From $I\cap\underline{C} = \emptyset$ we get $I\cap\underline{X} = I\cap\underline{Y} = \emptyset$. Also, by \cref{single_element_cocircuits}, $(X,\sigma(X))_j$ and $(Y,\sigma(Y))_j$ are cocircuits of $\mathcal{M}$. Hence if $i\notin\underline{X}$ or $i\notin\underline{Y}$, then $(I\setminus\{j\})\cup\{i\}$ is not a basis of $\mathcal{M}$ by \cref{defn_basis_rank}. Otherwise we have $X_i = Y_i \neq 0$ (since $X$ and $Y$ are conformal), and $X\neq Y$ (since $\sigma(X)\neq\sigma(Y)$). Then by (C3) of \cref{defn_cocircuits}, there exists a cocircuit of $\mathcal{M}$ whose support is contained in $(\underline{X}\cup\underline{Y}\cup\{j\})\setminus\{i\}\subseteq E\setminus((I\setminus\{j\})\cup\{i\})$, whence $(I\setminus\{j\})\cup\{i\}$ is not a basis of $\mathcal{M}$.
\end{pf}
\begin{cor}[geometric interpretation of $i\to_\epsilon j$-perturbation]\label{perturbation_interpretation}
Suppose that $V\in\Gr_{k,n}$, $i,j\in [n]$, and $\epsilon\in\{+,-\}$. For $\alpha\in\mathbb{R}$, let $W(\alpha)\in\Gr_{k,n}$ be the row span of the $k\times n$ matrix $[x^{(1)} | \cdots | x^{(j-1)} | (x^{(j)} + \alpha x^{(i)}) | x^{(j+1)} | \cdots | x^{(n)}]$, where $[x^{(1)} | \cdots | x^{(n)}]$ is a $k\times n$ matrix whose rows span $V$. (Note that $W(\alpha)$ does not depend on the choice of matrix.) Then for all $\alpha\in\mathbb{R}$ with sign $\epsilon$ such that $\Delta_I(W(\alpha))$ has the same sign as $\Delta_I(V)$ for all $I\in\binom{[n]}{k}$ with $\Delta_I(V)\neq 0$, $\mathcal{M}(W(\alpha))$ is the $i\to_\epsilon j$-perturbation of $\mathcal{M}(V)$.
\end{cor}
Note that the possible values of $\alpha$ form an open interval between $0$ and some number, or $\pm\infty$, with sign $\epsilon$.
\begin{eg}\label{perturbation_example}
Let $V\in\Gr_{2,4}$ be the row span of the matrix $\begin{bmatrix}1 & 0 & 2 & 0 \\ 0 & 3 & -1 & 4\end{bmatrix}$, and for $\alpha < 0$ let $W(\alpha)\in\Gr_{2,4}$ be the row span of the matrix $\begin{bmatrix}1 & 0 & 2 & \alpha \\ 0 & 3 & -1 & 4\end{bmatrix}$. Note that the $\{3,4\}$-minor of the first matrix equals $8$, and the $\{3,4\}$-minor of the second matrix equals $8 + \alpha$, so we should pick $\alpha > -8$ so that these minors agree in sign. In fact, for all $\alpha\in(-8,0)$ the corresponding minors of the two matrices agree in sign whenever the first minor is nonzero, whence $\mathcal{M}(W(\alpha))$ equals the $1\to_-4$-perturbation of $\mathcal{M}(V)$.
\end{eg}
\begin{pf}[of \cref{perturbation_interpretation}]
Note that for $I\in\binom{[n]}{k}$ and $\alpha\in\mathbb{R}$, we have
\begin{align}\label{perturbation_plueckers}
\Delta_I(W(\alpha)) = \begin{cases}
\Delta_I(V) + (-1)^{|I\cap (i,j)|}\alpha\Delta_{(I\setminus\{j\})\cup\{i\}}(V), & \text{if $i\notin I$ and $j\in I$} \\
\Delta_I(V), & \text{otherwise}
\end{cases},
\end{align}
where $(i,j)$ denotes the set of elements of $[n]$ strictly between $i$ and $j$. Hence the result follows from \cref{perturbation_chirotope}.
\end{pf}

We observe that certain $i\to_\epsilon j$-perturbations do not increase sign variation.
\begin{lem}[sign variation and $i\to_\epsilon j$-perturbation]\label{perturbation_sign_variation}
Suppose that $\mathcal{M}$ is an oriented matroid of rank $k$ with ground set $[n]$, and $m\ge k-1$. \\
(i) Let $\mathcal{N}$ be either the $(i+1)\to_+ i$-perturbation of $\mathcal{M}$ ($i\in [n-1]$), the $i\to_+ (i+1)$-perturbation of $\mathcal{M}$ ($i\in [n-1]$), the $1\to_{(-1)^m}n$-perturbation of $\mathcal{M}$, or the $n\to_{(-1)^m}1$-perturbation of $\mathcal{M}$. If $\var(X)\le m$ for all $X\in\mathcal{V}^*(\mathcal{M})$, then $\var(Y)\le m$ for all $Y\in\mathcal{V}^*(\mathcal{N})$. \\
(ii) Suppose that $\mathcal{P}\ge\mathcal{M}$ has rank $k$. If $\overline{\var}(X)\le m$ for all $X\in\mathcal{V}^*(\mathcal{M})\setminus\{0\}$, then $\overline{\var}(Y)\le m$ for all $Y\in\mathcal{V}^*(\mathcal{P})\setminus\{0\}$.
\end{lem}
Note that (i) above does not hold for any other $i\to_\epsilon j$-perturbations of $\mathcal{M}$ (assuming $i\neq j$); for counterexamples, we can take $k := 1$ and $m\in \{0,1\}$.
\begin{pf}
(i) Note that for $X\in\{0,+,-\}^n$, we have $\var(X) = \var((X_n, X_{n-1}, \cd, X_1))$, and if $\var(X)\le m$ then $\var((X_2, X_3, \cd, X_n, (-1)^mX_1))\le m$. By this cyclic symmetry, it will suffice to prove the result assuming that $\mathcal{N}$ is the $2\to_+ 1$-perturbation of $\mathcal{M}$. Suppose that $\var(X)\le m$ for all $X\in\mathcal{V}^*(\mathcal{M})$, but there exists a covector $Y$ of $\mathcal{N}$ with $\var(Y)\ge m+1$. We will derive a contradiction by showing that $Y$ is a covector of $\mathcal{M}$.

Since $\mathcal{M}\setminus\{1\} = \mathcal{N}\setminus\{1\}$, by \cref{defn_restriction} we have $X|_{[n]\setminus\{1\}} = Y|_{[n]\setminus\{1\}}$ for some covector $X$ of $\mathcal{M}$. From $\var(X)\le m$, we get that $Y_1\neq 0, Y_2$. Write $\mathcal{M}$ as the single element extension of $\mathcal{M}\setminus\{1\}$ by $\sigma:\mathcal{C}^*(\mathcal{M}\setminus\{1\})\to\{0,+,-\}$. Since $Y$ is a composition of conformal cocircuits of $\mathcal{N}$ (\cref{conformality_for_covectors}), by \cref{single_element_cocircuits} we have a composition of conformal cocircuits
$$
Y = (C^{(1)},(\sigma\circ\phi_2)(C^{(1)}))_1 \circ \cdots \circ (C^{(r)},(\sigma\circ\phi_2)(C^{(r)}))_1 \circ (D^{(1)}\circ E^{(1)},0)_1 \circ \cdots \circ (D^{(s)}\circ E^{(s)},0)_1
$$
for some cocircuits $C^{(1)}, \cd, C^{(r)}, D^{(1)}, E^{(1)}, \cd, D^{(s)}, E^{(s)}$ of $\mathcal{M}\setminus\{1\}$. If $Y_2 = 0$, then $(\sigma\circ\phi_2)(C^{(t)}) = \sigma(C^{(t)})$ for $t\in [r]$, whence by \cref{single_element_cocircuits}
\begin{multline*}
Y = (C^{(1)},\sigma(C^{(1)}))_1 \circ \cdots \circ (C^{(r)},\sigma(C^{(r)}))_1 \circ \\ (D^{(1)},\sigma(D^{(1)}))_1 \circ (E^{(1)},\sigma(E^{(1)}))_1\circ \cdots \circ (D^{(s)},\sigma(D^{(s)}))_1 \circ (E^{(s)},\sigma(E^{(s)}))_1
\end{multline*}
is a covector of $\mathcal{M}$, a contradiction. Hence $Y_1 = -Y_2$. In particular, because the composition above is of conformal cocircuits, we have $C^{(t)}_2\neq Y_1$ and $(\sigma\circ\phi_2)(C^{(t)})\neq -Y_1$ for $t\in [r]$. We also have $\sigma(C^{(u)}) = Y_1$ for some $u\in [r]$, since otherwise $(\sigma\circ\phi_2)(C^{(t)}) = \phi_2(C^{(t)}) \neq Y_1$ for all $t\in [r]$. This gives
\begin{multline*}
Y = (C^{(u)},\sigma(C^{(u)}))_1 \circ (C^{(1)},\sigma(C^{(1)}))_1 \circ \cdots \circ (C^{(r)},\sigma(C^{(r)}))_1 \circ \\ (D^{(1)},\sigma(D^{(1)}))_1 \circ (E^{(1)},\sigma(E^{(1)}))_1\circ \cdots \circ (D^{(s)},\sigma(D^{(s)}))_1 \circ (E^{(s)},\sigma(E^{(s)}))_1,
\end{multline*}
so $Y$ is a covector of $\mathcal{M}$ by \cref{single_element_cocircuits}, a contradiction.

(ii) This follows from a general fact about oriented matroids $\mathcal{A}$ and $\mathcal{B}$ with the same rank and ground set (7.7.5 of \cite{bjorner_las_vergnas_sturmfels_white_ziegler}): $\mathcal{A}\le\mathcal{B}$ iff for all nonzero covectors $Y$ of $\mathcal{B}$, there exists a nonzero covector $X$ of $\mathcal{A}$ with $X\le Y$.
\end{pf}

We now explain how to perturb an oriented matroid into a uniform oriented matroid by repeatedly applying $i\to_\epsilon j$-perturbations.
\begin{prop}[uniform perturbation]\label{uniform_perturbation}
Suppose that $\mathcal{M}$ is an oriented matroid of rank $k$ with ground set $[n]$. \\
(i) The oriented matroid obtained from $\mathcal{M}$ by applying any $k(2n-k-1)$ consecutive perturbations of the sequence
$$
\cd, (n-1)\to n, n\to 1, 1\to 2, 2\to 3, \cd, (n-1)\to n, n\to 1, 1\to 2, \cdots
$$
in order from left to right is uniform. (Here an $i\to j$-perturbation denotes either of the $i\to_\epsilon j$-perturbations, for $\epsilon\in \{+,-\}$.) \\
(ii) The oriented matroid obtained from $\mathcal{M}$ by applying any $(n-k)(n+k-1)$ consecutive perturbations of the sequence
$$
\cd, n\to (n-1), 1\to n, 2\to 1, 3\to 2, \cd, n\to (n-1), 1\to n, 2\to 1, \cdots
$$
in order from left to right is uniform. \\
(iii) The oriented matroid obtained from $\mathcal{M}$ by applying the sequence of perturbations
$$
1\to 2, 2\to 3, \cd, (n-1)\to n, n\to (n-1), (n-1)\to (n-2), \cd, 2\to 1
$$
in order from left to right $k$ times is uniform. \\
(iv) The oriented matroid obtained from $\mathcal{M}$ by applying the sequence of perturbations
$$
2\to 1, 3\to 2, \cd, n\to (n-1), (n-1)\to n, (n-2)\to (n-1), \cd, 1\to 2
$$
in order from left to right $n-k$ times is uniform.
\end{prop}
Thus we have four specific algorithms for perturbing $\mathcal{M}$ into a uniform oriented matroid, each using at most $2n^2$ perturbations. For example, if $k := 1$ and $n := 3$, then applying any of the following sequences of perturbations to $\mathcal{M}$, in order from left to right, produces a uniform oriented matroid: \\
$\bullet$ $1\to 2, 2\to 3, 3\to 1, 1\to 2$ (by (i)); or \\
$\bullet$ $2\to 3, 3\to 1, 1\to 2, 2\to 3$ (by (i)); or \\
$\bullet$ $3\to 1, 1\to 2, 2\to 3, 3\to 1$ (by (i)); or \\
$\bullet$ $2\to 1, 3\to 2, 1\to 3, 2\to 1, 3\to 2, 1\to 3$ (by (ii)); or \\
$\bullet$ $3\to 2, 1\to 3, 2\to 1, 3\to 2, 1\to 3, 2\to 1$ (by (ii)); or \\
$\bullet$ $1\to 3, 2\to 1, 3\to 2, 1\to 3, 2\to 1, 3\to 2$ (by (ii)); or \\
$\bullet$ $1\to 2, 2\to 3, 3\to 2, 2\to 1$ (by (iii)); or \\
$\bullet$ $2\to 1, 3\to 2, 2\to 3, 1\to 2, 2\to 1, 3\to 2, 2\to 3, 1\to 2$ (by (iv)).
\begin{eg}\label{eg_uniform_perturbation}
Let $V\in\Gr_{2,3}$ be the row span of the matrix $\begin{bmatrix}1 & 3 & 0 \\ 0 & 0 & 1\end{bmatrix}$, so that the vectors in $V$ change sign at most $m := 1$ time. Now $V$ is not generic, because $\Delta_{\{2,3\}}(V) = 0$. We can perturb $V$ into a generic subspace by applying a $3\to_- 1$-perturbation, giving the row span of
$$
\begin{bmatrix}
1 & 3 & 0 \\
\alpha & 0 & 1
\end{bmatrix} \qquad (\alpha < 0),
$$
or by applying a $3\to_+2$-perturbation, giving the row span of
$$
\begin{bmatrix}
1 & 3 & 0 \\
0 & \beta & 1
\end{bmatrix} \qquad (\beta > 0).
$$
The vectors in either of these generic subspaces change sign at most once, as guaranteed by \cref{perturbation_sign_variation}. Note that we cannot make $V$ generic by applying only $1\to 2$- and $2\to 3$-perturbations.
\end{eg}

\begin{pf}[of \cref{uniform_perturbation}]
Let $\mathcal{N}$ be an oriented matroid of rank $k$ with ground set $[n]$. The {\itshape dual} $\mathcal{N}^*$ of $\mathcal{N}$ is the oriented matroid of rank $n-k$ with ground set $[n]$ whose chirotope is given by $\chi_{\mathcal{N}^*}(J) = (-1)^{\sum_{j\in J}j}\chi_\mathcal{N}([n]\setminus J)$ for $J\in\binom{[n]}{n-k}$. Note that $\mathcal{N}$ is uniform iff $\mathcal{N}^*$ is uniform. Also, \cref{perturbation_chirotope} implies that the dual of the $i\to_\epsilon j$-perturbation of $\mathcal{N}$ is the $j\to_{-\epsilon}i$-perturbation of $\mathcal{N}^*$. Hence statements (i) and (ii) are dual, and statements (iii) and (iv) are dual. We will prove (ii) and (iv).

A {\itshape hyperplane} of an (oriented) matroid is a maximal subset of its ground set which contains no basis. Note that by \cref{defn_basis_rank} and (C2) of \cref{defn_cocircuits}, hyperplanes are precisely the complements of supports of cocircuits. Now suppose that we have a collection of functions, each of which, given an oriented matroid $\mathcal{P}$ of rank $k$ with ground set $[n]$, produces an oriented matroid $\mathcal{P}'\ge\mathcal{P}$ of rank $k$, such that no hyperplane of $\mathcal{P}$ of size at least $k$ is a hyperplane of $\mathcal{P}'$. Note that every basis of $\mathcal{P}$ is a basis of $\mathcal{P}'$ (by \cref{defn_poset}), so every hyperplane of $\mathcal{P}'$ is contained in a hyperplane of $\mathcal{P}$. Hence the maximum size of a hyperplane of $\mathcal{P}'$ is less than the maximum size of a hyperplane of $\mathcal{P}$, unless every hyperplane of $\mathcal{P}$ has size less than $k$ (i.e.\ $\mathcal{P}$ is uniform). By applying such a function $n-k$ times (possibly a different function in our collection each time), we obtain a uniform oriented matroid. Thus to prove (ii), it suffices to show that for all $i\in\mathbb{Z}$, applying the sequence of perturbations
$$
(i+1) \to i, (i+2)\to (i+1), \cd, (i+n+k-1)\to (i+n+k-2)
$$
(where we read the indices modulo $n$) in order from left to right is such a function (we then apply this function for $i = j, j+(n+k-1), j+2(n+k-1), \cd, j+(n-k-1)(n+k-1)$ for any $j\in\mathbb{Z}$). Similarly, to prove (iv), it suffices to show that applying the sequence of perturbations
$$
2\to 1, 3\to 2, \cd, n\to (n-1), (n-1)\to n, (n-2)\to (n-1), \cd, 1\to 2
$$
in order from left to right is such a function. To this end, we prove the following claim.
\begin{claim}
Suppose that $\mathcal{P}$ is an oriented matroid of rank $k$ with ground set $[n]$, and $I\subseteq [n]$ is a hyperplane of $\mathcal{P}$ with $|I| \ge k$. Take $a\in I$ and $b\in [k]$ such that \\
$\bullet$ $a$ is not a coloop of $\mathcal{P}|_I$; \\
$\bullet$ $a+1, a+2, \cd, a+b-1\in I$ are coloops of $\mathcal{P}|_I$; and \\
$\bullet$ $a+b\notin I$, \\
where we read the indices modulo $n$. Then for all $\mathcal{Q}\ge\mathcal{P}$ of rank $k$, $I$ is not a hyperplane the oriented matroid obtained from $\mathcal{Q}$ by applying the sequence of perturbations $(a+1)\to a, (a+2)\to(a+1), \cd, (a+b)\to(a+b-1)$ in order from left to right, where we read the indices modulo $n$.
\end{claim}

\begin{claimpf}
First note that $I\neq [n]$, and $\mathcal{P}|_I$ has at most $k-1$ coloops (otherwise the rank of $\mathcal{P}|_I$ would be at least $k$), so such $a$ and $b$ exist. Also note that for any oriented matroids $\mathcal{A}\le\mathcal{B}$ of equal rank with ground set $[n]$, by \cref{perturbation_chirotope} the $i\to_\epsilon j$-perturbation of $\mathcal{A}$ is less than or equal to the $i\to_\epsilon j$-perturbation of $\mathcal{B}$, for all $i,j\in [n]$ and $\epsilon\in\{+,-\}$. Hence it will suffice to prove the claim assuming that $\mathcal{Q} = \mathcal{P}$.

Let $\mathcal{P}^{(0)} := \mathcal{P}$, and define $\mathcal{P}^{(c)}$ recursively for $c = 1, \cd, b$ as either of the $(a+c)\to_\epsilon (a+c-1)$-perturbations of $\mathcal{P}^{(c-1)}$ for $\epsilon\in\{+,-\}$. Also let $J\in\binom{I}{k-1}$ be a basis of $\mathcal{P}|_I$ which does not contain $a$. Since $a+1, a+2, \cd, a+b-1$ are coloops of $\mathcal{P}|_I$, they are in $J$.

We claim that $(J\cup\{a,a+b\})\setminus\{a+c\}$ is a basis of $\mathcal{P}^{(c)}$ for $0 \le c \le b$. Let us prove this by induction on $c$. For the base case $c = 0$, we must show that $J\cup\{a+b\}$ is a basis of $\mathcal{P}$. If not, then by \cref{defn_basis_rank} there exists a cocircuit $C$ of $\mathcal{P}$ with $(J\cup\{a+b\})\cap\underline{C} = \emptyset$. We then have $C_I = 0$. (Otherwise there exists a cocircuit $D$ of $\mathcal{P}_I$ with $D\le C_I$ by \cref{conformality_for_covectors}, whence $J\cap\underline{D} = \emptyset$, contradicting \cref{defn_basis_rank} since $J$ is a basis of $\mathcal{P}_I$.) This gives $\underline{C}\subset [n]\setminus I$, which contradicts (C2) of \cref{defn_cocircuits} because $[n]\setminus I$ is the support of a cocircuit of $\mathcal{P}$. For the induction step, suppose that $c\in [b]$ and $(J\cup\{a,a+b\})\setminus\{a+c-1\}$ is a basis of $\mathcal{P}^{(c-1)}$. By \cref{perturbation_chirotope},  we have
\begin{multline*}
\chi_{\mathcal{P}^{(c)}}((J\cup\{a,a+b\})\setminus\{a+c\}) = \\
\begin{cases}
\pm\chi_{\mathcal{P}^{(c-1)}}((J\cup\{a,a+b\})\setminus\{a+c-1\}), & \text{if $\chi_{\mathcal{P}^{(c-1)}}((J\cup\{a,a+b\})\setminus\{a+c\}) = 0$} \\
\chi_{\mathcal{P}^{(c-1)}}((J\cup\{a,a+b\})\setminus\{a+c\}), & \text{otherwise}
\end{cases}.
\end{multline*}
In the first case we have $\chi_{\mathcal{P}^{(c)}}((J\cup\{a,a+b\})\setminus\{a+c\})\neq 0$ by the induction hypothesis, while in the second case $(J\cup\{a,a+b\})\setminus\{a+c\}$ is a basis of $\mathcal{P}^{(c-1)}$, and hence also of $\mathcal{P}^{(c)}\ge\mathcal{P}^{(c-1)}$. This completes the induction. Taking $c := b$ we get that $J\cup\{a\}$ is a basis of $\mathcal{P}^{(c)}$, and so $I$ is not a hyperplane of $\mathcal{P}^{(c)}$.
\end{claimpf}

Note that for any $a\in\mathbb{Z}$ and $b\in [k]$, the sequence $(a+1)\to a, (a+2)\to (a+1), \cd, (a+b)\to (a+b-1)$ is a consecutive subsequence of
$$
(i+1)\to i, (i+2)\to (i+1), \cd, (i+n+k-1)\to (i+n+k-2)
$$
for all $i\in\mathbb{Z}$ (where we read the indices modulo $n$). This proves (ii).

For (iv), let $\mathcal{P}$ be an oriented matroid of rank $k$ with ground set $[n]$, and $I\subseteq [n]$ a hyperplane of $\mathcal{P}$ with $|I|\ge k$. It will suffice to show that $I$ is not a hyperplane of the oriented matroid $\mathcal{P}'$ obtained from $\mathcal{P}$ by applying the sequence of perturbations
$$
2\to 1, 3\to 2, \cd, n\to (n-1), (n-1)\to n, (n-2)\to (n-1), \cd, 1\to 2
$$
in order from left to right. To this end, take $i\in [n]\setminus I$. If there exists an element of $[1,i]\cap I$ which is not a coloop of $\mathcal{P}|_I$, then we may take $a$ and $b$ as in the statement of the claim such that we also have $1 \le a < a+b \le i$; then $I$ is not a hyperplane of the oriented matroid obtained from any $\mathcal{Q}\ge\mathcal{P}$ by applying the sequence of perturbations $(a+1)\to a, (a+2)\to (a+1), \cd, (a+b)\to (a+b-1)$ in order from left to right, whence $I$ is not a hyperplane of $\mathcal{P}'$. Otherwise, there exists an element of $[i,n]\cap I$ which is not a coloop of $\mathcal{P}|_I$, whence we take $a'\in [i,n]\cap I$ and $b'\in [k]$ such that \\
$\bullet$ $a'$ is not a coloop of $\mathcal{P}|_I$; \\
$\bullet$ $a'-1, a'-2, \cd, a'-b'+1\in I$ are coloops of $\mathcal{P}|_I$; and \\
$\bullet$ $a'-b'\notin I$. \\
We have $a'-b' \ge i$, and by the claim $I$ is not a hyperplane of the oriented matroid obtained from any $\mathcal{Q}\ge\mathcal{P}$ by applying the sequence of perturbations $(a'-1)\to a', (a'-2)\to (a'-1), \cd, (a'-b')\to (a'-b'+1)$ in order from left to right, whence $I$ is not a hyperplane of $\mathcal{P}'$.
\end{pf}

We are now ready to give a necessary and sufficient condition that $\var(X)\le m$ for all $X\in\mathcal{V}^*(\mathcal{M})$.
\begin{thm}\label{perturbation_criterion}
Suppose that $\mathcal{M}$ is an oriented matroid of rank $k$ with ground set $[n]$, and $m\ge k-1$. Let $\mathcal{N}$ be any oriented matroid obtained from $\mathcal{M}$ by applying one of the following sequences of perturbations: \\
$\bullet$ any $k(2n-k-1)$ consecutive perturbations of the sequence
$$
\cd, (n-1)\to_+ n, n\to_{(-1)^m} 1, 1\to_+ 2, 2\to_+ 3, \cd, (n-1)\to_+ n, n\to_{(-1)^m} 1, 1\to_+ 2, \cdots
$$
in order from left to right; or \\
$\bullet$ any $(n-k)(n+k-1)$ consecutive perturbations of the sequence
$$
\cd, n\to_+ (n-1), 1\to_{(-1)^m}n, 2\to_+ 1, 3\to_+ 2, \cd, n\to_+ (n-1), 1\to_{(-1)^m}n, 2\to_+ 1, \cdots
$$
in order from left to right; or \\
$\bullet$ the sequence of perturbations
$$
1\to_+ 2, 2\to_+ 3, \cd, (n-1)\to_+ n, n\to_+ (n-1), (n-1)\to_+ (n-2), \cd, 2\to_+ 1
$$
in order from left to right $k$ times; or \\
$\bullet$ the sequence of perturbations
$$
2\to_+ 1, 3\to_+ 2, \cd, n\to_+ (n-1), (n-1)\to_+ n, (n-2)\to_+ (n-1), \cd, 1\to_+ 2
$$
in order from left to right $n-k$ times.

Then $\mathcal{N}$ is uniform, and the following are equivalent: \\
(i) $\var(X)\le m$ for all $X\in\mathcal{V}^*(\mathcal{M})$; \\
(ii) $\var(Y)\le m$ for all $Y\in\mathcal{V}^*(\mathcal{N})$; and \\
(iii) $\var((\chi_\mathcal{N}(I\cup\{i\}))_{i\in [n]\setminus I})\le m-k+1$ for all $I\in\binom{[n]}{k-1}$.
\end{thm}
Note that the first two sequences of perturbations take advantage of the cyclic symmetry of sign variation, but they depend on (the parity of) $m$, whereas the last two sequences do not. Note that none of the sequences depend on $\mathcal{M}$ (only on $n$ and $k$, and perhaps $m$).
\begin{pf}
\cref{uniform_perturbation} implies that $\mathcal{N}$ is uniform. We have (i) $\Rightarrow$ (ii) by \cref{perturbation_sign_variation}, (ii) $\Rightarrow$ (i) by \cref{defn_poset}, (ii) $\Rightarrow$ (iii) by \cref{chirotope_criterion}(i), and (iii) $\Rightarrow$ (ii) by \cref{chirotope_criterion}(ii) (since $\mathcal{N}$ is uniform).
\end{pf}
We can interpret this statement as a closure result in the space of oriented matroids (or the Grassmannian $\Gr_{k,n}$), where the {\itshape closure} of a set $S$ of oriented matroids is $\{\mathcal{M} : \mathcal{M}\le\mathcal{N} \text{ for some } \mathcal{N}\in S\}$. ($\Gr_{k,n}$ has the classical topology.)
\begin{thm}\label{generic_elements_are_dense}
Let $n\ge k \ge 0$ and $m\ge k-1$. \\
(i) Let $S$ be the set of oriented matroids $\mathcal{M}$ of rank $k$ with ground set $[n]$ satisfying $\var(X)\le m$ for all $X\in\mathcal{V}^*(\mathcal{M})$. Then the closure of the set of uniform elements of $S$ (in the space of oriented matroids of rank $k$ with ground set $[n]$) equals $S$. \\
(ii) Let $T := \{V\in\Gr_{k,n} : \var(v)\le m \text{ for all } v\in V\}$. Then the closure in $\Gr_{k,n}$ of the set of generic elements of $T$ equals $T$.
\end{thm}

\begin{pf}
\cref{perturbation_criterion} implies (i). For (ii), note that the closure in $\Gr_{k,n}$ of the generic elements of $T$ is contained in $T$. Conversely, given $V\in T$ we can construct (by \cref{perturbation_criterion}) a sequence $\mathcal{M}_0 := \mathcal{M}(V), \mathcal{M}_1, \mathcal{M}_2, \cd, \mathcal{M}_r$ of elements of $S$ such that $\mathcal{M}_s$ is the $i_s\to_{\epsilon_s} j_s$-perturbation of $\mathcal{M}_{s-1}$ (for some $i_s,j_s\in [n]$ and $\epsilon_s\in\{+,-\}$) for all $s\in [r]$, and $\mathcal{M}_r$ is uniform. For $\alpha > 0$, let $V_0(\alpha) := V$, and define $V_s(\alpha)\in\Gr_{k,n}$ recursively for $s = 1, \cd, r$ as the row span of the $k\times n$ matrix $[x^{(1)} | \cdots | x^{(j_s-1)} | (x^{(j_s)} + \epsilon_s\alpha^{2^{s-1}}x^{(i_s)}) | x^{(j_s+1)} | \cdots | x^{(n)}]$, where $[x^{(1)} | \cdots | x^{(n)}]$ is a $k\times n$ matrix whose rows span $V_{s-1}(\alpha)$. Note that for $0 \le s \le r$, every Pl\"{u}cker coordinate of $V_s(\alpha)$ is a polynomial in $\alpha$ of degree at most $2^s-1$; we can prove this by induction on $s$, using \rf{perturbation_plueckers}.
\begin{claim}
Let $s\in [r]$ and $I\in\binom{[n]}{k}$. Then for $\alpha > 0$ sufficiently small, either $\Delta_I(V_{s-1}(\alpha)) = 0$, or $\Delta_I(V_s(\alpha))$ and $\Delta_I(V_{s-1}(\alpha))$ are nonzero with the same sign.
\end{claim}

\begin{claimpf}
Regard $\Delta_I(V_{s-1}(\alpha))$ as a polynomial in $\alpha$. If this polynomial is zero then the claim is proven, so suppose that this polynomial is nonzero, and write $\Delta_I(V_{s-1}(\alpha)) = c\alpha^d + O(\alpha^{d+1})$ (as $\alpha\to 0$) for some $d \le 2^{s-1} - 1$ and $c \neq 0$. Then by \rf{perturbation_plueckers} we have $\Delta_I(V_s(\alpha)) = \Delta_I(V_{s-1}(\alpha)) + O(\alpha^{2^{s-1}}) = c\alpha^d + O(\alpha^{d+1})$. Hence for $\alpha > 0$ sufficiently small, we have $\sign(\Delta_I(V_s(\alpha))) = \sign(\Delta_I(V_{s-1}(\alpha))) = \sign(c)$.
\end{claimpf}
Thus by \cref{perturbation_interpretation}, for $\alpha > 0$ sufficiently small we have $\mathcal{M}(V_s(\alpha)) = \mathcal{M}_s$ for all $s\in [r]$, whence $V_r(\alpha)$ is generic and $V_r(\alpha)\in T$. Taking $\alpha\to 0$ shows explicitly that $V$ is in the closure of $T$.
\end{pf}

\section{Defining amplituhedra and Grassmann polytopes}\label{sec_amplituhedron}

\noindent Let $k,n,r\in\mathbb{N}$ with $n\ge k,r$, and let $Z:\mathbb{R}^n\to\mathbb{R}^r$ be a linear map, which we also regard as an $r\times n$ matrix. Arkani-Hamed and Trnka \cite{arkani-hamed_trnka} consider the map $\Gr_{k,n}^{\ge 0}\to\Gr_{k,r}$ induced by $Z$ on the totally nonnegative Grassmannian. Explicitly, if $X$ is a $k\times n$ matrix whose row span is $V\in\Gr_{k,n}^{\ge 0}$, then $Z(V)$ is the row span of the $k\times r$ matrix $XZ^T$. In the case that $k\le r$ and all $r\times r$ minors of $Z$ are positive, Arkani-Hamed and Trnka call the image of this map a {\itshape (tree) amplituhedron}, and use it to calculate scattering amplitudes in $\mathcal{N} = 4$ supersymmetric Yang-Mills theory (taking $r := k+4$). One motivation they provide for requiring that $k\le r$ and $Z$ have positive $r\times r$ minors is to guarantee that the map $\Gr_{k,n}^{\ge 0}\to\Gr_{k,r}$ induced by $Z$ is well defined, i.e.\ that $Z(V)$ has dimension $k$ for all $V\in\Gr_{k,n}^{\ge 0}$. As a more general sufficient condition for this map to be well defined, Lam \cite{lam} requires that the row span of $Z$ has a $k$-dimensional subspace which is totally positive. (It is not obvious that Arkani-Hamed and Trnka's condition is indeed a special case of Lam's; see Section 15.1 of \cite{lam}.) In the case that the map $\Gr_{k,n}^{\ge 0}\to\Gr_{k,r}$ induced by $Z$ is well defined, Lam calls the image a {\itshape (full) Grassmann polytope}, since in the case $k = 1$ Grassmann polytopes are precisely polytopes in the projective space $\Gr_{1,r} = \mathbb{P}^{r-1}$ (and the amplituhedra are projective {\itshape cyclic polytopes}). In this section we give (\cref{amplituhedron_map}) a necessary and sufficient condition for the map $\Gr_{k,n}^{\ge 0}\to\Gr_{k,r}$ to be well defined, in terms of sign variation; we are able to translate this into a condition on the maximal minors of $Z$ using the results of \cref{sec_sign_changes}. As a consequence, we recover Arkani-Hamed and Trnka's and Lam's sufficient conditions. To be thorough, we similarly determine when the map $\Gr_{k,n}^{>0}\to\Gr_{k,r}$ induced by $Z$ on the totally positive Grassmannian is well defined (\cref{amplituhedron_map_2}).

\begin{lem}\label{vector_extension}
Let $v\in\mathbb{R}^n\setminus\{0\}$ and $k\le n$. \\
(i) There exists an element of $\Gr_{k,n}^{\ge 0}$ containing $v$ iff $\var(v)\le k-1$. \\
(ii) There exists an element of $\Gr_{k,n}^{>0}$ containing $v$ iff $\overline{\var}(v)\le k-1$.
\end{lem}

\begin{pf}
The forward directions of (i) and (ii) follow from Gantmakher and Krein's result (\cref{gantmakher_krein}). For the reverse direction of (i), suppose that $\var(v)\le k-1$. Then we may partition $[n]$ into pairwise disjoint nonempty intervals of integers $I_1, \cd, I_k$, such that for all $j\in [k]$ the components of $v|_{I_j}$ are all nonnegative or all nonpositive. For $j\in [k]$, let $w^{(j)}\in\mathbb{R}^n$ have support $I_j$ such that $w^{(j)}|_{I_j}$ equals $v|_{I_j}$ if $v|_{I_j}\neq 0$, and $e_{I_j}$ otherwise. Then $\spn(\{w^{(j)} : j\in [k]\})\in\Gr_{k,n}^{\ge 0}$ contains $v$. (For example, if $v = (2,5,0,-1,-4,-1,0,0,3)$ and $k = 4$, then we may take $I_1 := \{1,2,3\}, I_2 := \{4,5,6\}, I_3 := \{7,8\}, I_4 := \{9\}$, whence our subspace is the row span of the matrix
$$
\begin{bmatrix}
2 & 5 & 0 & 0 & 0 & 0 & 0 & 0 & 0 \\
0 & 0 & 0 & -1 & -4 & -1 & 0 & 0 & 0 \\
0 & 0 & 0 & 0 & 0 & 0 & 1 & 1 & 0 \\
0 & 0 & 0 & 0 & 0 & 0 & 0 & 0 & 3
\end{bmatrix};
$$
note that $v$ is the sum of rows $1$, $2$, and $4$.)

Now we prove the reverse direction of (ii). The point is that by rescaling the basis vectors of $\mathbb{R}^n$ (the {\itshape torus action} on the Grassmannian), we need only determine the sign vectors appearing in totally positive subspaces. 
\begin{claim}[\cite{gantmakher_krein_1950}, \cite{bjorner_las_vergnas_sturmfels_white_ziegler}] Let $V\in\Gr_{k,n}^{>0}$. \\
(i) $\{\sign(v) : v\in V\} = \{X\in\{0,+,-\}^n : \overline{\var}(X) \le k-1\}\cup\{0\}$. \\
(ii) $\{\sign(w) : w\in V^\perp\} = \{X\in\{0,+,-\}^n : \var(X)\ge k\}\cup\{0\}$.
\end{claim}
In the terminology of oriented matroids, the sets in (i) and (ii) are the covectors and vectors, respectively, of $\mathcal{M}(V)$.
\begin{claimpf}
This essentially follows from known results, as follows. First recall that by \cref{dual_translation}(ii), $V$ is totally positive iff $\alt(V^\perp)$ is totally positive. Hence by \cref{dual_translation}(i), parts (i) and (ii) of the claim are equivalent. Let us prove (ii). The containment $\subseteq$ follows from Gantmakher and Krein's result (\cref{dual_gantmakher_krein}(ii)). For the containment $\supseteq$, given $X\in\{0,+,-\}^n$ with $\var(X)\ge k$, take $I\in\binom{[n]}{k+1}$ such that $X$ alternates in sign on $I$. By Proposition 9.4.1 of \cite{bjorner_las_vergnas_sturmfels_white_ziegler}, there exists $w\in V^\perp$ such that $\sign(w|_I) = X|_I$ and $\sign(w|_{[n]\setminus I}) = 0$. Now for each $j\in [n]\setminus I$, take $v^{(j)}\in V^\perp$ such that $v^{(j)}_j = 1$ and $v^{(j)}|_{[n]\setminus (I\cup \{j\})} = 0$. (For example, fix any $h\in I$, whence $\Delta_{([n]\setminus I)\cup\{h\}}(V^\perp)\neq 0$ since $V^\perp$ is generic. Then take any $(n-k)\times n$ matrix whose rows span $V^\perp$, and row reduce it so that we get an identity matrix in the columns $([n]\setminus I)\cup\{h\}$. Then we let $v^{(j)}$ for $j\in [n]\setminus I$ be the row of this matrix whose pivot column is $j$.) By perturbing $w$ by $v^{(j)}$ for $j\in [n]\setminus I$ so that $w_j$ has sign $X_j$, we obtain a vector in $V^\perp$ with sign vector $X$.
\end{claimpf}

Suppose that $\overline{\var}(v)\le k-1$. Take any $V\in\Gr_{k,n}^{>0}$ (e.g.\ let $V$ be the row span of the matrix
$$
\begin{bmatrix}
1 & 1 & \cdots & 1 \\
t_1 & t_2 & \cdots & t_n \\
t_1^2 & t_2^2 & \cdots & t_n^2 \\
\vdots & \vdots & \ddots & \vdots \\
t_1^{k-1} & t_2^{k-1} & \cdots & t_n^{k-1}
\end{bmatrix},
$$
where $t_1 < \cdots < t_n$). Then the oriented matroid $\mathcal{M}(V)$ defined by $V$ is the alternating oriented matroid of rank $k$ with ground set $[n]$, whence $\sign(v)$ is a covector of $\mathcal{M}$ by the claim. That is (cf.\ \cref{defn_realizable_oriented_matroids}), there exist $\alpha_1, \cd, \alpha_n > 0$ such that $(\alpha_1v_1, \cd, \alpha_nv_n)\in V$. Then $\{(\frac{w_1}{\alpha_1}, \cd, \frac{w_n}{\alpha_n}) : w\in V\}\in\Gr_{k,n}^{>0}$ contains $v$.
\end{pf}

\begin{thm}\label{amplituhedron_map}
Suppose that $k,n,r\in\mathbb{N}$ with $n\ge k,r$, and that $Z:\mathbb{R}^n\to\mathbb{R}^{r}$ is a linear map, which we also regard as an $r\times n$ matrix. Let $d$ be the rank of $Z$ and $W\in\Gr_{d,n}$ the row span of $Z$, so that $W^\perp = \ker(Z)\in\Gr_{n-d,n}$. The following are equivalent: \\
(i) the map $\Gr_{k,n}^{\ge 0}\to\Gr_{k,r}$ induced by $Z$ is well defined, i.e.\ $\dim(Z(V)) = k$ for all $V\in\Gr_{k,n}^{\ge 0}$; \\
(ii) $\var(v)\ge k$ for all nonzero $v\in\ker(Z)$; and \\
(iii) $\overline{\var}((\Delta_{I\setminus\{i\}}(W))_{i\in I})\le d-k$ for all $I\in\binom{[n]}{d+1}$ such that $W|_I$ has dimension $d$.
\end{thm}
We explain how to use \cref{amplituhedron_map} to deduce the sufficient conditions of Arkani-Hamed and Trnka, and of Lam, for the map $\Gr_{k,n}^{\ge 0}\to\Gr_{k,r}$ induced by $Z$ to be well defined. Note that if the $r\times r$ minors of $Z$ are all positive, then $d = r$ and $W$ is totally positive, so the condition (iii) holds for any $k\le r$. Alternatively, by \cref{dual_gantmakher_krein}(ii), we have $\var(v)\ge r$ for all nonzero $v\in\ker(Z)$, so the condition (ii) holds for any $k\le r$. This recovers the sufficient condition of Arkani-Hamed and Trnka \cite{arkani-hamed_trnka}. On the other hand, if $W$ has a subspace $V\in\Gr_{k,n}^{>0}$, then by \cref{dual_gantmakher_krein}(ii) we have $\var(v)\ge k$ for all $v\in V^\perp\setminus\{0\}$, which implies condition (ii) above since $\ker(Z) = W^\perp\subseteq V^\perp$. This recovers the sufficient condition of Lam \cite{lam}. However, our result does not show why Arkani-Hamed and Trnka's condition is a special case of Lam's. Indeed, it is an interesting open problem to determine whether or not Lam's sufficient condition is also necessary, i.e.\ whether the condition $\var(v)\ge k$ for all nonzero $v\in W^\perp$ implies that $W$ has a totally positive $k$-dimensional subspace.

\begin{eg}
Let $Z:\mathbb{R}^4\to\mathbb{R}^2$ be the linear map given by the matrix $\begin{bmatrix}2 & -1 & 1 & 1 \\ 1 & 2 & -1 & 3\end{bmatrix}$ (so $n = 4$, $d = r = 2$), and let $W\in\Gr_{2,4}$ be the row span of this matrix. Let us use \cref{amplituhedron_map}(iii) to determine for which $k$ ($0 \le k \le 4$) the map $\Gr_{k,4}^{\ge 0}\to\Gr_{k,2}$ induced by $Z$ is well defined. The $4$ relevant sequences of Pl\"{u}cker coordinates (as $I$ ranges over $\binom{[4]}{3}$) are
\begin{align*}
& (\Delta_{\{2,3\}}(W), \Delta_{\{1,3\}}(W), \Delta_{\{1,2\}}(W)) = (-1, -3, 5), \\
& (\Delta_{\{2,4\}}(W), \Delta_{\{1,4\}}(W), \Delta_{\{1,2\}}(W)) = (-5, 5, 5), \\
& (\Delta_{\{3,4\}}(W), \Delta_{\{1,4\}}(W), \Delta_{\{1,3\}}(W)) = (4, 5, -3), \\
& (\Delta_{\{3,4\}}(W), \Delta_{\{2,4\}}(W), \Delta_{\{2,3\}}(W)) = (4, -5, -1).
\end{align*}
The maximum number of sign changes among these $4$ sequences is $1$, which is at most $2-k$ iff $k\le 1$. Hence the map is well defined iff $k\le 1$.

Note that for $k\ge 2$, the proof of \cref{vector_extension}(i) shows how to explicitly construct $V\in\Gr_{k,4}^{\ge 0}$ with $\dim(Z(V)) < k$: take a nonzero $v\in\ker(Z)$ with $\var(v)\le 1$, and extend $v$ to $V\in\Gr_{k,4}^{\ge 0}$. For example, if $k = 2$ we can take $v = (1,-3,-5,0)\in\ker(Z)$ and extend it to the row span $V\in\Gr_{2,4}^{\ge 0}$ of the matrix $\begin{bmatrix}1 & 0 & 0 & 0 \\ 0 & -3 & -5 & 0\end{bmatrix}$. Note that $Z(V)$ is the span of $(2,1)$, so $\dim(Z(V)) = 1 < \dim(V)$.
\end{eg}

\begin{pf}[of \cref{amplituhedron_map}]
(i) $\Leftrightarrow$ (ii): The map $\Gr_{k,n}^{\ge 0}\to\Gr_{k,r}$ induced by $Z$ is well defined iff for all $V\in\Gr_{k,n}^{\ge 0}$ and $v\in V\setminus\{0\}$, we have $Z(v)\neq 0$. This condition is equivalent to (ii) above by \cref{vector_extension}(i).

(ii) $\Leftrightarrow$ (iii): This is precisely the dual statement of (the realizable case of) \cref{chirotope_criterion}(ii). Explicitly, recall that $\alt:\mathbb{R}^n\to\mathbb{R}^n$ is defined by $\alt(v) := (v_1, -v_2, v_3, \cd, (-1)^{n-1}v_n)$ for $v\in\mathbb{R}^n$. By \cref{dual_translation}(i), the condition (ii) is equivalent to $\overline{\var}(v)\le n-k-1$ for all $v\in\alt(\ker(Z))\setminus\{0\}$, which is in turn equivalent to $\overline{\var}((\Delta_{J\cup\{i\}}(\alt(\ker(Z))))_{i\in [n]\setminus J})\le d-k$ for all $J\in\binom{[n]}{n-d-1}$ such that $\Delta_{J\cup\{i\}}(\alt(\ker(Z)))\neq 0$ for some $i\in [n]$ by \cref{chirotope_criterion}(ii). This condition is precisely (iii) above, since $\Delta_K(W) = \Delta_{[n]\setminus K}(\alt(\ker(Z)))$ for all $K\in\binom{[n]}{d}$ by \cref{dual_translation}(ii).
\end{pf}

We give the analogue of \cref{amplituhedron_map} for the map induced by $Z$ not on $\Gr_{k,n}^{\ge 0}$, but on $\Gr_{k,n}^{>0}$.
\begin{thm}\label{amplituhedron_map_2}
Suppose that $k,n,r\in\mathbb{N}$ with $n\ge k,r$, and that $Z:\mathbb{R}^n\to\mathbb{R}^{r}$ is a linear map, which we also regard as an $r\times n$ matrix. Let $d$ be the rank of $Z$ and $W\in\Gr_{d,n}$ the row span of $Z$, so that $W^\perp = \ker(Z)\in\Gr_{n-d,n}$. The following are equivalent: \\
(i) the map $\Gr_{k,n}^{>0}\to\Gr_{k,r}$ induced by $Z$ is well defined, i.e.\ $\dim(Z(V)) = k$ for all $V\in\Gr_{k,n}^{>0}$; \\
(ii) $\overline{\var}(v)\ge k$ for all nonzero $v\in\ker(Z)$; and \\
(iii) there exists a generic perturbation $W'\in\Gr_{d,n}$ of $W$ such that $\var((\Delta_{I\setminus\{i\}}(W'))_{i\in I})\le d-k$ for all $I\in\binom{[n]}{d+1}$.
\end{thm}
We omit the proof, since it is similar to that of \cref{amplituhedron_map}; we only mention that instead of \cref{vector_extension}(i) we use \cref{vector_extension}(ii), and along with \cref{chirotope_criterion}(ii) we also use \cref{generic_elements_are_dense}.

\section{Positroids from sign vectors}\label{sec_positroids}

\noindent Recall that the totally nonnegative Grassmannian $\Gr_{k,n}^{\ge 0}$ has a cell decomposition, where the {\itshape positroid cell} of $V\in\Gr_{k,n}^{\ge 0}$ is determined by $M(V) := \{I\in\binom{[n]}{k} : \Delta_I(V)\neq 0\}$. The goal of this section is show how to obtain the positroid cell of a given $V\in\Gr_{k,n}^{\ge 0}$ from the sign vectors of $V$ (i.e.\ $\mathcal{V}^*(\mathcal{M}(V))$). Note that $M(V)$ is the set of bases of $\mathcal{M}(V)$, so $\mathcal{V}^*(\mathcal{M}(V))$ determines $M(V)$ by the theory of oriented matroids. However, this does not exploit the fact that $V$ is totally nonnegative. We now describe two other ways to recover $M(V)$ from the sign vectors of $V$, both of which require $V$ to be totally nonnegative.

We begin by examining the {\itshape Schubert cell} of $V$, which is labeled by the lexicographic minimum of $M(V)$. Recall that the {\itshape Gale partial order} $\le_\textnormal{Gale}$ on $\binom{[n]}{k}$ is defined by
$$
I\le_\textnormal{Gale}J \quad\iff\quad i_1 \le j_1,\; i_2 \le j_2,\; \cd,\; i_k \le j_k
$$
for subsets $I = \{i_1, \cd, i_k\}$ ($i_1 < \cdots < i_k$), $J = \{j_1, \cd, j_k\}$ ($j_1 < \cdots < j_k$) of $[n]$. Note that $I\le_\textnormal{Gale}J$ iff $|I\cap [m]|\ge |J\cap [m]|$ for all $m\in [n]$. Also recall that for $V\in\Gr_{k,n}$, $A(V)$ is the set of $I\in\binom{[n]}{k}$ such that some vector in $V$ strictly alternates in sign on $I$. Note that if $I\in M(V)$ then $V|_I = \mathbb{R}^I$, so $M(V)\subseteq A(V)$. We can obtain the Schubert cell of $V\in\Gr_{k,n}^{\ge 0}$ from $A(V)$ as follows.
\begin{thm}[Schubert cell from sign vectors]\label{schubert_cell_criterion}
For $V\in\Gr_{k,n}^{\ge 0}$, the lexicographic minimum of $M(V)$ equals the Gale minimum of $A(V)$.
\end{thm}
We remark that the lexicographic minimum of $M(V)$ is also the Gale minimum of $M(V)$, for all $V\in\Gr_{k,n}$. (In general, the lexicographically minimal basis of any matroid with a totally ordered ground set is also a Gale minimum \cite{gale}.) However, $A(V)$ does not necessarily equal $M(V)$ (see \cref{alternating_set_can_be_bigger} or \cref{same_alternating_set}), nor does $A(V)$ necessarily uniquely determine $M(V)$ (see \cref{same_alternating_set}). Also, if $V$ is not totally nonnegative, then $A(V)$ does not necessarily have a Gale minimum (see \cref{gale_minimality_example}).
\begin{eg}\label{same_alternating_set}
Let $V,W\in\Gr_{2,3}^{\ge 0}$ be the row spans of the matrices $\begin{bmatrix}1 & 0 & -1 \\ 0 & 1 & 0\end{bmatrix}$, $\begin{bmatrix}1 & 0 & -1 \\ 0 & 1 & 1\end{bmatrix}$, respectively. Then $A(V) = A(W) = \binom{[3]}{2}$, but $M(V)\neq M(W)$ since $\{1,3\}\in M(W)\setminus M(V)$.
\end{eg}
\begin{eg}\label{gale_minimality_example}
Let $V\in\Gr_{3,6}$ be the row span of the matrix $\begin{bmatrix}1 & 0 & -1 & -1 & 1 & 0 \\ 0 & 1 & 1 & 2 & 0 & 0 \\ 0 & 0 & 0 & 0 & 0 & 1\end{bmatrix}$, which is not totally nonnegative. Then $(1,-1,-2,-3,1,0)\in V$ strictly alternates in sign on $\{1,2,5\}$, and $(3,2,-1,1,3,0)\in V$ strictly alternates in sign on $\{1,3,4\}$, but no vector in $V$ strictly alternates in sign on $\{1,2,3\}$ or $\{1,2,4\}$. Hence $A(V)$ has no Gale minimum.
\end{eg}
\begin{pf}
Given $V\in\Gr_{k,n}^{\ge 0}$, let $I$ be the lexicographic minimum of $M(V)$.
\begin{claim}
Let $m\in [n]$ and $l := |I\cap [m]|$. Then $V|_{[m]}\in\Gr_{l,m}^{\ge 0}$.
\end{claim}

\begin{claimpf}
Express $V$ as the row span of a $k\times n$ matrix $X = [x^{(1)} | \cdots | x^{(n)}]$ whose restriction to the columns in $I$ is an identity matrix. Note that $V|_{[m]}$ is the row span of the first $m$ columns of $X$. Since $\{x^{(i)} : i\in I\cap [m]\}$ is linearly independent, we may extend $I\cap [m]$ to $B\in M(V|_{[m]})$, and then extend $B$ to $B'\in M(V)$. Since $I$ is the Gale minimum of $M(V)$ \cite{gale}, we have $I\le_\textnormal{Gale}B'$. In particular $|I\cap [m]|\ge |B'\cap [m]|$, so $B = I\cap [m]$. Hence $\dim(V|_{[m]}) = l$, and the entries in the first $m$ columns of $X$ past the $l$th row are all zero. It follows that $V|_{[m]}$ is the row span of the submatrix of $X$ formed by the first $l$ rows and the first $m$ columns. Since the restriction of $X$ to the columns in $I$ is an identity matrix, we see that $\Delta_K(V|_{[m]}) = \Delta_{K\cup(I\setminus [m])}(V)\ge 0$ for $K\in\binom{[m]}{l}$.
\end{claimpf}
Hence if $v\in V$ strictly alternates in sign on $J\in\binom{[n]}{k}$, by \cref{gantmakher_krein} we get $|I\cap [m]| - 1\ge \var(v|_{[m]})\ge |J\cap [m]| - 1$ for all $m\in [n]$, whence $I\le_\textnormal{Gale}J$.
\end{pf}

Given $n\ge 0$, for $j\in [n]$ let $\le_j$ be the total order on $[n]$ defined by $j <_j j+1 <_j \cdots <_j n <_j 1 <_j \cdots <_j j-1$. Then for $V\in\Gr_{k,n}$, we let $I_j$ ($j\in [n]$) denote the lexicographic minimum of $M(V)$ with respect to $\le_j$. The tuple $(I_1, \cd, I_n)$ is called the {\itshape Grassmann necklace} of $V$. For example, if $V\in\Gr_{2,4}$ is generic, then the Grassmann necklace of $V$ is $(\{1,2\}, \{2,3\}, \{3,4\}, \{4,1\})$. The Grassmann necklace is of special interest to us because of a result of Postnikov (Theorem 17.1 of \cite{postnikov}), which implies that if $V$ is totally nonnegative, the positroid cell of $V$ is determined by its Grassmann necklace. Oh \cite{oh} explicitly described $M(V)$ in terms of the Grassmann necklace of $V$, for $V\in\Gr_{k,n}^{\ge 0}$. 
\begin{thm}[\cite{oh}]\label{oh}
Suppose that $V\in\Gr_{k,n}^{\ge 0}$ has Grassmann necklace $(I_1, \cd, I_n)\in\binom{[n]}{k}^n$. Then
\begin{align*}
M(V) = \left\{J\in\binom{[n]}{k} : I_j \le_\textnormal{$j$-Gale}J \textnormal{ for all } j\in [n]\right\}.
\end{align*}

\end{thm}
(Here $\le_\textnormal{$j$-Gale}$ denotes the Gale order on $\binom{[n]}{k}$ induced by $\le_j$.)

We can generalize \cref{schubert_cell_criterion} to the Grassmann necklace $(I_1, \cd, I_n)$ of $V\in\Gr_{k,n}^{\ge 0}$ as follows. For $j\in [n]$, we define $V_j$ as the row span of the cyclically shifted $k\times n$ matrix $[x^{(j)} | x^{(j+1)} | \cdots | x^{(n)} | (-1)^{k-1}x^{(1)} | \cdots | (-1)^{k-1}x^{(j-1)}]$, where $[x^{(1)} | \cdots | x^{(n)}]$ is a $k\times n$ matrix whose rows span $V$. Note that $V_j$ does not depend on our choice of matrix, and since $V$ is totally nonnegative so is $V_j$. Then $\{i-j+1\pmod{n} : i\in I_j\}$ is the lexicographic minimum of $M(V_j)$, and so applying \cref{schubert_cell_criterion} to $V_j$ gives the following result.
\begin{cor}[Grassmann necklace from sign vectors]\label{grassmann_necklace_criterion}
Suppose that $V\in\Gr_{k,n}^{\ge 0}$. For $j\in [n]$, let $A_j$ be the set of $J\in\binom{[n]}{k}$ such that some vector in $V$ strictly alternates in sign on $J$ except precisely from component $\max(J\cap[1,j))$ to component $\min(J\cap[j,n])$ (if both components exist). Then $A_j$ has a $j$-Gale minimum $I_j$ for all $j\in [n]$, and $(I_1, \cd, I_n)$ is the Grassmann necklace of $V$.
\end{cor}
For example, if $n := 5$, $J := \{1,3,4,5\}$, and $j := 3$, then $(1,1,1,-1,1)$ strictly alternates in sign on $J$ except precisely from component $\max(J\cap[1,j))$ to component $\min(J\cap[j,n])$, but $(1,1,-1,1,-1)$ does not. (If $j\le \min(J)$ or $j > \max(J)$, then the condition reduces to ``strictly alternates in sign on $J$.'')

With Oh's result (\cref{oh}), we get the following corollary.
\begin{cor}\label{basis_criterion}
Suppose that $V\in\Gr_{k,n}^{\ge 0}$ has Grassmann necklace $(I_1, \cd, I_n)\in\binom{[n]}{k}^n$, and $J\in\binom{[n]}{k}$. Then the following are equivalent: \\
(i) $J\in M(V)$; \\
(ii) $I_j\le_\textnormal{$j$-Gale}J$ for all $j\in [n]$; and \\
(iii) $V$ realizes all $2k$ sign vectors in $\{+,-\}^J$ which alternate in sign between every pair of consecutive components, with at most one exceptional pair.
\end{cor}
For example, if $k = 5$ the $2k$ sign vectors in (iii) above are $(+,-,+,-,+)$, $(+,+,-,+,-)$, $(+,-,-,+,-)$, $(+,-,+,+,-)$, $(+,-,+,-,-)$, and their negations. Since $V$ realizes a sign vector iff $V$ realizes its negation, we need only check $k$ sign vectors in (iii) up to sign.
\begin{pf}
We have (i) $\Rightarrow$ (iii) since $J\in M(V)$ implies $V|_J = \mathbb{R}^J$, (iii) $\Rightarrow$ (ii) by \cref{grassmann_necklace_criterion}(ii), and (ii) $\Rightarrow$ (i) by Oh's result (\cref{oh}).
\end{pf}
We can prove (iii) $\Rightarrow$ (i) directly from \cref{gantmakher_krein}, as follows. Suppose that (iii) holds, but $J\notin M(V)$. Then there exists $v\in V\setminus\{0\}$ with $v|_J = 0$; take $j\in [n]$ such that $v_j\neq 0$. Then (iii) guarantees the existence of a vector $w\in V$ which strictly alternates in sign on $J$ except precisely from component $\max(J\cap[1,j))$ to component $\min(J\cap[j,n])$ (if both components exist). Adding a sufficiently large multiple of $\pm v$ to $w$ gives a vector in $V$ which strictly alternates in sign on $J\cup\{j\}$, contradicting \cref{gantmakher_krein}. This establishes the equivalence of (i) and (iii) without appleaing to Oh's result (\cref{oh}). The implication (i) $\Rightarrow$ (ii) is a general fact about matroids \cite{gale}. We would be interested to see a direct proof of (ii) $\Rightarrow$ (iii) (and hence of \cref{basis_criterion}) which is substantially different from Oh's proof, using the tools of sign variation.

\begin{rmk}\label{check_all_patterns}
We remark that (iii) $\Rightarrow$ (i) does not necessarily hold when $V$ is not totally nonnegative; in fact, it is possible that $V$ realizes all $2^k$ sign vectors in $\{+,-\}^J$ except two, but $J\notin M(V)$. To see this, given $J\in\binom{[n]}{k}$, let $v\in\mathbb{R}^J$ have no zero components, and take $V\in\Gr_{k,n}$ such that $V|_J = \{v\}^\perp$ (which is always possible, assuming $n > k$). That is, $J\notin M(V)$ and $V|_J = \{w\in\mathbb{R}^J : \sum_{j\in J}v_jw_j = 0\}$. We see that if $w\in\mathbb{R}^J$ satisfies $\sign(w) = \sign(v)$, then $\sum_{j\in J}v_jw_j > 0$, and so $w\notin V|_J$. Similarly, if $\sign(w) = -\sign(v)$ then $w\notin V|_J$. Conversely, given $\omega\in\{+,-\}^J$ with $\omega\neq\pm\sign(v)$, let us construct $w\in V|_J$ with $\sign(w) = \omega$. Take $a,b\in J$ such that $\sign(v_a)\omega_a \neq \sign(v_b)\omega_b$. For $j\in J\setminus\{a,b\}$ let $w_j$ be any real number with sign $\omega_j$, then take $w_b$ with sign $\omega_b$ and sufficiently large magnitude that $\sign(\sum_{j\in J\setminus\{a\}}v_jw_j) = \sign(v_b)\omega_b$, and set $w_a := -\frac{\sum_{j\in J\setminus\{a\}}v_jw_j}{v_a}$. Thus $V$ realizes all sign vectors in $\{+,-\}^J$ except for precisely $\pm\sign(v)$.

On the other hand, if $V$ realizes all $2^k$ sign vectors in $\{+,-\}^J$, then $J\in M(V)$. Indeed, if $J\notin M(V)$ then we may take $v\in (V|_J)^\perp\setminus\{0\}$, whence $V$ does not realize any $\omega\in\{+,-\}^J$ satisfying $\sign(v)\le\omega$.
\end{rmk}

\bibliographystyle{math}
\bibliography{Covectors_bib}

\end{document}